\documentclass[12pt]{amsart}
\usepackage[utf8]{inputenc}
\usepackage{hyperref}
\usepackage{mathtools}
\usepackage[dvipsnames]{xcolor}
\usepackage{amsmath, amssymb, amsthm, tikz, float}
\usepackage{array}
\usepackage{arydshln}
\usepackage{comment}
\usepackage{multirow}
\usepackage{geometry}
\usepackage{lineno}
\usepackage{cancel}
\theoremstyle{plain}
\usetikzlibrary{decorations.markings}
\newtheorem{theorem}{Theorem}[section] 
\newtheorem{lemma}[theorem]{Lemma}
\newtheorem{proposition}[theorem]{Proposition}
\newtheorem{corollary}[theorem]{Corollary}
\newtheorem{remark}[theorem]{Remark}
\theoremstyle{definition}
\newtheorem{example}[theorem]{Example}
\newtheorem{definition}[theorem]{Definition}
\newtheorem{identity}[theorem]{Identity}

\newcommand{\bmbar}[1]{\mbox{\boldmath $\overline #1$}}

\DeclareMathOperator{\Der}{Der}  
\DeclareMathOperator{\Flip}{Flip}  
\DeclareMathOperator{\R}{\mathcal{R}}  
\DeclareMathOperator{\Fib}{{Fib}}  
\DeclareMathOperator{\Cat}{{Cat}}  
\DeclareMathOperator{\Jac}{{Jac}}  
\DeclareMathOperator{\Pas}{{Pas}}  
\DeclareMathOperator{\Sha}{{Sha}}  

\definecolor{darkblue}{rgb}{0.0, 0.0, 0.7}
\definecolor{darkgreen}{rgb}{0.0, 0.5, 0.2}

\title[On Sums, Derivatives, and Flips of Riordan Arrays]{On Sums, Derivatives, and Flips of Riordan Arrays}
\author[Bang, von Bell, Culver, Dickson, Dimitrov, Perrier, Sundaram]{Caroline Bang \and Matias von Bell \and Eric Culver \and Jessica Dickson \and Stoyan Dimitrov \and Rachel Perrier \and Sheila Sundaram }

\address{Caroline Bang: Iowa State University, 411 Morrill Rd, Ames, IA 50011.}
\email{cbang@iastate.edu}
\address{Matias von Bell: 
Institute of Geometry, Graz  University of Technology, Kopernikusgasse 24, Graz, A-8010, Austria.}
\email{matias.vonbell@gmail.com}
\address{Eric Culver: University of Colorado Denver, 1201 Larimer St., Denver, CO 80204}
\email{eric.culver@ucdenver.edu}
\address{Jessica Dickson: Washington State University, PO Box 643113, Pullman, WA 99164, USA}
\email{jmdickson@wsu.edu}
\address{Stoyan Dimitrov: Rutgers University, Hill Center, Piscataway, NJ 08854, USA}
\email{EmailToStoyan@gmail.com}
\address{Rachel Perrier: Washington State University, PO Box 643113, Pullman, WA 99164, USA}
\email{rachel.perrier@wsu.edu}
\address{Sheila Sundaram: Pierrepont School, One Sylvan Road North, Westport, CT 06880, USA}
\email{shsund@comcast.net}
\date{\today}

\begin{document}


\maketitle
\begin{abstract}
We study three operations on Riordan arrays.
First, we investigate when the sum of Riordan arrays yields another Riordan array.
We characterize the $A$- and $Z$-sequences of these sums of Riordan arrays, and also identify an analog for $A$-sequences when the sum of Riordan arrays does not yield a Riordan array.
In addition, we define the new operations `$\Der$' and `$\Flip$' on Riordan arrays.  
We fully characterize the Riordan arrays resulting from these operations applied to the Appell and Lagrange subgroups of the Riordan group.
Finally, we study the application of these operations to various known Riordan arrays, generating many combinatorial identities in the process. 
\end{abstract}
\tableofcontents

\section{Introduction}

Riordan arrays are lower triangular matrices, extending infinitely to the right and downward, whose columns encode generating functions. They were first introduced in \cite{shapiro1991riordan} (see also \cite{Shapiro2005}) as a generalization of Pascal’s triangle, which satisfies the matrix equation:

$$
\begin{bmatrix}
1 & 0 & 0 & 0 & 0 & \cdots \\
1 & 1 & 0 & 0 & 0 & \cdots \\
1 & 2 & 1 & 0 & 0 & \cdots \\
1 & 3 & 3 & 1 & 0 & \cdots \\
1 & 4 & 6 & 4 & 1 & \cdots \\
\vdots & \vdots & \vdots & \vdots & \vdots & \ddots\\
\end{bmatrix}
\begin{bmatrix}
1 \\ x \\ x^2 \\ x^3 \\ x^4 \\ \vdots 
\end{bmatrix}
=
\begin{bmatrix}
1 \\ 1+x \\ (1+x)^2 \\ (1+x)^3 \\ (1+x)^4 \\ \vdots 
\end{bmatrix}
$$

In general, a Riordan array $\R(d(t),h(t))$ is characterized by a pair of formal power series $d(t)$ and $h(t)$, satisfying $d(0)\neq 0$ and $h'(0)\neq 0$.
Riordan arrays form a group under matrix multiplication known as the Riordan group, which has been extensively studied \cite{Barry, shapiro2003bijections, shapiro1991riordan}. 
However,  Riordan arrays are not in general closed under matrix addition. 
In Section~\ref{sec:Sums} we investigate the sums of Riordan arrays, and characterize when the sum of two Riordan arrays is a Riordan array.

\vspace{0.3cm}
\noindent \textbf{Theorem~\ref{thm:R+S}}
{\em
Let $R = \mathcal{R}(d_R(t), h_R(t))$ and $S = \mathcal{R}(d_S(t), h_S(t))$ be two Riordan arrays. Then $R+S$ is a Riordan array if and only if $h_R(t)=h_S(t)$ and $d_R(0) + d_S(0) \ne 0$. In this case $R+S$ is the Riordan array $\mathcal{R}(d_R(t) + d_S(t), h(t))$.
}
\vspace{0.3cm}

Instead of using power series, Riordan arrays can also be described in terms of a pair of sequences $A$ and $Z$. 
In the case that Riordan arrays can be added, we characterize the $A$- and $Z$-sequences of their sums (Theorem~\ref{ZSeqR+S}). 
When two Riordan arrays do not have the same second argument  $h(t)$, their sum is not a Riordan array. 
However, we can still consider the array formed by their matrix sum $R+S$. 
We call such arrays \textit{Riordan sumrays} and show that they satisfy a second order recurrence (Theorem~\ref{thm:Eric}).

In Section~\ref{sec:derflip}, we introduce two new operations on Riordan arrays which we call {\em derivatives} and {\em flips}. 
For a Riordan array $\R(d(t),h(t))$, derivative and flip are the following two operations respectively.
\begin{enumerate}
    \item $\Der$: $\R(d(t), h(t)) \mapsto \R(h'(t), t d(t));$
   \item $\Flip$: $\R(d(t), h(t)) \mapsto \R(h(t)/t, t d(t))$.  
\end{enumerate}

An important observation is that $\Flip$ is an involution, and that $\Der$ maps the Appell subgroup of the Riordan group to the Lagrange subgroup, and vice versa. 
In the case of Appell and Lagrange subgroups, we can fully characterize the derivatives as follows.

\vspace{0.3cm}
\noindent \textbf{Theorem~\ref{thm:derAppell}}
{\em
For a Riordan array $\R(d(t),t)$ in the Appell subgroup and $m\geq 0$ we have 
$$ \Der^{2m}(R) = \R\left(\sum_{i= 0}^m S_{m+1,i+1} t^i d^{(i)}(t),t\right) = \R\left(\sum_{i\geq 0} d_i(i+1)^mt^i, t \right),$$
$$ \Der^{2m+1}(R) = \R\left(1, \sum_{i= 0}^m S_{m+1,i+1} t^{i+1} d^{(i)}(t)\right) = \R\left(1, t\sum_{i\geq 0} d_i(i+1)^mt^i \right),$$
where $S_{m+1,i+1}$ is the Stirling number of the second kind, and $d^{(i)}(t)$ denotes the $i$th derivative of $d(t)$.   Here $\Der^k$ (respectively, $\Flip^k$)  denotes the operation $\Der$ (respectively, $\Flip$) iterated $k$ times. 
}

There is a large literature on Riordan arrays showing their use in  obtaining combinatorial identities \cite{LuzonMerliniMoronSprugnoli2012LAA,shapiro1991riordan,SprugnoliFundThm1994,sprugnoli2011combinatorial}. 
In the spirit of the pioneering work of Louis Shapiro, we  highlight this feature of Riordan arrays throughout the article. The new operations $\Der$ and $\Flip$ have  proved especially fruitful by revealing surprising connections with known sequences, as well as leading to the discovery of new identities.
These identities arise primarily from viewing a Riordan array $\R$ as the transformation of a sequence $(a_i)_{i=0}^\infty$ to a sequence $(b_i)_{i=0}^\infty$ via the product $\R\cdot (a_i)_{i=0}^\infty = (b_i)_{i=0}^\infty$. 
In Section~\ref{sec:applications}, we apply our results about $\Der$ and $\Flip$ to various Riordan arrays including the Fibonacci array, Pascal array, Catalan array, and Shapiro array. 
These are respectively the Riordan arrays
\begin{align*} \Fib &\coloneqq \R(1, t+ t^2),
&\Pas &\coloneqq  \R\left(\dfrac{1}{1-t}, \dfrac{t}{1-t}\right), \\ 
\Cat &\coloneqq  \R(C(t), tC(t)),  
&\Sha &\coloneqq  \R(C(t), tC^2(t)),
\end{align*}
where $C(t)$ is the generating function for the Catalan numbers (see Table~\ref{table:sequences}). 
As a consequence we obtain several known combinatorial identities in new ways, along with identities that we have not found in the literature  (Identities \ref{id:2}, \ref{id:riordanMotzkin}, \ref{id:powerOf4}, \ref{id:OddsToCat}, \ref{id:pureDescents}, \ref{id:avoiding}, \ref{id:weakCompCatalan2}, \ref{id:derSH}, and Theorem~\ref{thm:latticePaths}). 
For several of the identities obtained we also provide new combinatorial proofs (Identities \ref{id:2},\ref{id:riordanMotzkin}, \ref{id:OddsToCat}, \ref{id:pureDescents}, \ref{id:derSH}). 
Furthermore, for any Riordan array in the Appell subgroup, the operation $\Der$ can be interpreted as giving a Riordan array involving weighted integer compositions (see Theorem~\ref{thm:DerAsComps}).
A summary of the relationships between sequences obtained using various Riordan arrays is given in Table~\ref{table:SummaryTable}.

\begin{table}
\centering
\def\arraystretch{2.2}
    \begin{tabular}{|l|l|l|}
        \hline
        Name & Sequence & Generating Function 
        \\
        \hline
        Fibonacci & $(F_i)_{i=0}^\infty = (1,1,2,3,5,8,13,\ldots)$ & $F(t) = \dfrac{1}{1-t-t^2}$ \rule[-0.4cm]{0pt}{0pt}\\
        \hline
        Catalan & $(C_i)_{i=0}^\infty = (1,1,2,5,14,42,\ldots)$ & $C(t) = \dfrac{1- \sqrt{1-4t}}{2t}$ 
        \rule[-0.4cm]{0pt}{0pt}\\
        \hline
        Motzkin & $(M_i)_{i=0}^\infty = (1,1,2,4,9,21,\ldots)$ &  $M(t)=\dfrac{1-t-\sqrt{1-2t-3t^2}}{2t^2}$ 
        \rule[-0.4cm]{0pt}{0pt}\\
        \hline
        Riordan & $(R_i)_{i=0}^\infty = (1,0,1,1,3,6,15,\ldots)$ & $R(t)=\dfrac{1+t-\sqrt{1-2t-3t^2}}{2t(1+t)}$ 
        \rule[-0.4cm]{0pt}{0pt}\\
        \hline
    \end{tabular}
    \vspace{0.2cm}
    \caption{The recurring sequences in this article and their generating functions.}
    \label{table:sequences}
\end{table}

\begin{table}
\begin{center}
\def\arraystretch{1.2}
\begin{tabular}{|p{0.6cm}|l|l|l|}
\hline
& Riordan array $\R$ & Sequence $(a_i)_{i=0}^\infty$ & Product $\R\cdot \,(a_i)_{i=0}^\infty$ \\
\hline
1 & \multirow{2}{4cm}{$T:=\R(1+t,t)$}  & 
A000012, $(1)_{i=0}^\infty$ & A040000, $(1,2,2,\ldots)$ \\  
\cline{1-1} \cline{3-4}
2 & & A007598, $(F_i^2)_{i=0}^\infty$ & A001519, $(F_{2i})_{i=0}^\infty$  \\ 
\hline

3 & $\Fib :=\Der T = \R(1,t+t^2)$ & 
A000012, $(1)_{i=0}^\infty$ & A000045, $(F_i)_{i=0}^\infty$ \\  
\hline

4 
& $\Jac:= \Der^3 T = \R(1,t+2t^2)$ & 
A000012, $(1)_{i=0}^\infty$ & $($A001045$(i+1))_{i=0}^\infty$ \\  
\hline

5 
& \multirow{2}{3cm}{$\R\left(\frac{C(t)-1}{tF(t)},t\right)$}  & A000045, $(F_i)_{i=0}^\infty$ & A000108, $(C_i)_{i=1}^\infty$ \\

\cline{1-1} \cline{3-4}

6 
& & A000007, $(1,0,0,\ldots)$ & $(1,1, A067324(i)_{i=0}^\infty)$ \\
\hline

7 & \multirow{2}{3cm}{$\R(tM(t)+1,t)$}  & 
A000012, $(1)_{i=0}^\infty$ & $($A086615$(i)+1)_{i=-1}^\infty$ \\  
\cline{1-1} \cline{3-4}

8 
& & A005043, $(R_i)_{i=0}^\infty$ & A001006, $(M_i)_{i=0}^\infty$ \\
\hline

9 & \multirow{2}{5cm}{$R := \R\left(\frac{-1}{1-t^2}, \frac{t}{1-t^2}\right)$} & A000012, $(1)_{i=0}^\infty$ & A152163, $(-F_i)_{i=0}^\infty$ \\
\cline{1-1} \cline{3-4}
10 & & A000027, $(i)_{i=0}^\infty$ & $($A029907$(i))_{i=1}^\infty$ \\
\hline

11 & \multirow{2}{5cm}{$S := \R\left(\frac{2}{(t^2-1)^2}, \frac{t}{1-t^2}\right)$} & A000012, $(1)_{i=0}^\infty$ & $($A052952$(i))_{i=1}^\infty$ \\
\cline{1-1} \cline{3-4}
12 & & A000027, $(i+1)_{i=0}^\infty$ & $(2\cdot $A001629$(i))_{i=2}^\infty$ \\
\hline

13 & \multirow{2}{5cm}{$R+S$} & A000012, $(1)_{i=0}^\infty$ & $($A001350$(i))_{i=1}^\infty$ \\
\cline{1-1} \cline{3-4}
14 & & A000027, $(i+1)_{i=0}^\infty$ & A045925 \\
\hline

15 & $\Flip\Der(\Pas)$ & A000012, $(1)_{i=0}^\infty$ & A001519, $(F_{2i})_{i=0}^\infty$ \\
\hline

16 & $\Der^2(\Pas)$ & A000012, $(1)_{i=0}^\infty$ & A001906, $(F_{2i+1})_{i=0}^\infty$ \\
\hline

17 & $\Der^3(\Pas)$ & A000012, $(1)_{i=0}^\infty$ & A004146 \\
\hline

18 & $\Der^4(\Pas)$ & A000012, $(1)_{i=0}^\infty$ & A033453 \\
\hline

19 & $\Der^6(\Pas)$ & A000012, $(1)_{i=0}^\infty$ & A144109 \\
\hline

20 & $\Der^2(\Pas) + \Der^3(\Pas)$ & A000012, $(1)_{i=0}^\infty$ & A027941, $(F_{2i+1}-1)_{i=1}^\infty$ \\
\hline

21 & \multirow{3}{5cm}{$\Der(\Cat)$} & A000012, $(1)_{i=0}^\infty$ & A001700 \\
\cline{1-1} \cline{3-4}
22 & & A000027, $(2^i)_{i=0}^\infty$ & A000302, $(4^i)_{i=1}^\infty$ \\
\cline{1-1} \cline{3-4}
23 & & A005408, $(2i+1)_{i=0}^\infty$ & $($A129869$(i))_{i=1}^\infty$ \\
\hline

24 
& $\Flip\Der(\Cat)$ & A000012, $(1)_{i=0}^\infty$ & A026737 \\
\hline

25 & $\Der^2(\Cat)$ & A000012, $(1)_{i=0}^\infty$ & A026671 \\
\hline

26 & \multirow{2}{5cm}{$\Der(\Sha)$} & A000012, $(1)_{i=0}^\infty$ & A002054 \\
\cline{1-1} \cline{3-4}
27 
& & A000027, $(2^i)_{i=0}^\infty$ & A008549  \\
\hline
\end{tabular}
\end{center}
\vspace{.1in}
\caption{A summary of the relationships between sequences appearing in this article. 
}
\label{table:SummaryTable}
\end{table}


\section{Background and Examples}\label{sec:Background}

We begin by providing the necessary background on Riordan arrays.           
\begin{definition}
\label{def:riordan}
A {\em Riordan array} $D=\mathcal{R}(d(t),h(t))$ is an infinite, lower triangular array defined by a pair of formal power series $(d(t),h(t))$ where $d(0)\neq 0, h(0)=0,$ and $h'(0)\neq 0.$ The $(n,k)$-entry of $D$ is given by 
\begin{equation} \label{d}
d_{n,k} = [t^n]\{d(t)h(t)^k\}, n,k\ge 0,
\end{equation}
where $[t^n]\{a(t)\}=a_n$ is notation to extract the coefficient $a_n$ from an ordinary generating function $a(t)=\sum_{n\geq 0} a_n t^n$.  In particular note that $d_{0,0}=d(0)$ and $d_{0,k}=0$ for $k\ge 1$. 
\end{definition}

For example, choosing $d(t) = 1/(1-t)$ and $h(t)=t/(1-t)$ gives the {\em Pascal array} 
$$
\Pas=\mathcal{R} \left(\frac{1}{1-t}, \frac{t}{1-t}\right) = \begin{bmatrix}
1 & 0 & 0 & 0 & 0 & \cdots \\
1 & 1 & 0 & 0 & 0 & \cdots \\
1 & 2 & 1 & 0 & 0 & \cdots \\
1 & 3 & 3 & 1 & 0 & \cdots \\
1 & 4 & 6 & 4 & 1 & \cdots \\
\vdots & \vdots & \vdots & \vdots & \vdots & \ddots\\
\end{bmatrix}.
$$

We will encounter more examples in Section~\ref{subsec:examples}. 
First, however, we review the main theorems on Riordan arrays. 

\subsection{The main theorems}
\label{subsec:mainTheorems}
As mentioned in the introduction, Riordan arrays have been used in the literature to investigate combinatorial identities. 
These results often exploit the following theorem, often called the \emph{fundamental theorem of Riordan arrays} or ``FTRA" \cite{Barry}.

\begin{theorem}(FTRA) \label{thm:FTRA}\cite[Equation  (5)]{shapiro1991riordan}, \cite[Theorem 1.1]{SprugnoliFundThm1994} \\
Let $D=\mathcal{R}(d(t),h(t))$ be a Riordan array, and let $f(t):=\sum_{k\ge 0} f_k t^k$ be the generating function for the sequence  $(f_k)_{k\ge 0}.$ If $D$ maps the sequence $(f_k)_{k\ge 0},$ to the sequence $(g_k)_{k\ge 0},$ then the generating function $g(t)=\sum_{k\ge 0} g_k t^k$ for the sequence 
$(g_k)_{k\ge 0}$ is 
 \[d(t) f(h(t)).\]
\end{theorem}

We will provide several examples of using the FTRA to obtain combinatorial identities in Section~\ref{subsec:examples}.

An important fact about Riordan arrays is that they form a group known as the {\em Riordan group}. 
This is encapsulated in the following theorem.

\begin{theorem}\cite{shapiro1991riordan} 
\label{thm:Shapiro}
The Riordan arrays form a group with respect to matrix multiplication.  If $R_1 = \mathcal{R}(d_1(t), h_1(t))$ and $R_2 = \mathcal{R}(d_2(t), h_2(t))$ are Riordan arrays, then 
\[R_1R_2=\mathcal{R}(d_1(t) d_2(h_1(t)), h_2(h_1(t))).\]
The inverse of $R = \mathcal{R}(d(t), h(t))$  is
$$R^{-1}=\mathcal {R}\left(\dfrac{1}{d(\bmbar h(t))}, \bmbar h(t)\right),$$
where $\bmbar h(t)$ denotes the compositional inverse of $h(t),$ 
i.e., $h(\bmbar h (t))=t=\bmbar h (h(t))$.  The identity of the group is $\mathcal{R}(1, t)$.

\end{theorem}

The following are well-known subgroups of the Riordan group.

\begin{itemize}
    \item The {\em Lagrange subgroup}  $\{\R(1,h(t)) \mid h'(0)\neq 0\}$ (also called the Associated subgroup). 
    \item The {\em Appell subgroup} $\{\R(d(t),t) \mid d(0)\neq 0\}$ (also called the Toeplitz subgroup).
    \item The {\em Bell subgroup} $\{\R(d(t), t\cdot d(t)) \mid d(0) \neq 0\}$ (these are also called  Renewal arrays or Rogers arrays).
    \item The {\em derivative subgroup} $\{\R(h'(t), h(t)) \mid h'(0)\neq 0\}$.
    \item The {\em hitting time subgroup} $\{\R(th'(t)/h(t), h(t)) \mid h'(0)\neq 0\}$.
\end{itemize}

Associated to a Riordan array are two sequences, its $A$-sequence \cite{Rogers1978} and its $Z$-sequence \cite{Merlini-Rogers-Sprugnoli-Verri1997}.  The  $A$-sequence $A=(a_0 \neq 0,a_1,a_2,...)$ provides a way to recursively construct each row of the array from the previous row, via the generating function $A(t)=\sum_{n\geq0}a_nt^n$, as follows. 
\begin{equation*}\label{A}
d_{r+1,c+1}=a_0 d_{r,c}+a_1 d_{r,c+1}+a_2 d_{r,c+2}+ \cdots, \ r\ge 0.
\end{equation*} 

On the other hand, the $Z$-sequence, $Z=(z_0,z_1,z_2, ...)$, with generating function $Z(t)=\sum_{n\geq0}z_nt^n$, allows the $0$th column to be (uniquely) constructed using 
\begin{equation*} \label{B}
d_{r+1,0}=z_0 d_{r,0}+z_1 d_{r,1}+z_2 d_{r,2} + \cdots.    
\end{equation*} Note that since Riordan arrays are lower triangular,  
the sums in equations \eqref{A} and \eqref{B} are finite.
Results of \cite{Rogers1978},\cite[Theorem 2.1, Theorem 2.2]{HeSprugnoli2009} and \cite{Merlini-Rogers-Sprugnoli-Verri1997} imply that a Riordan array can be completely characterized by the triple ($d_{0,0},A(t), Z(t)$). In the case of the Riordan array $\Pas$, $A=(1,1,0,0,\dots)$ and thus  $A(t)=1+t,$ so that equation \eqref{A} reduces to the familiar Pascal recurrence $d_{r+1,c+1}= d_{r,c}+d_{r,c+1}$. Similarly, $Z=(1,0,0,\dots)$ and $Z(t)=1=d_{0,0}$ gives, using equation \eqref{B}, $d_{r+1,0}=1$.

\begin{theorem}(\cite[Theorem 2.1]{HeSprugnoli2009}, \cite[Theorem 2.2]{Merlini-Rogers-Sprugnoli-Verri1997}, \cite{Rogers1978})  
\label{thm:He-Sprugnoli}

\noindent If $R = \mathcal{R}(d(t), h(t))$ is a Riordan array, then the $A$-sequence is determined by:
\begin{equation*}
    \label{Aseq} h(t)=t A(h(t)) 
 \iff A(t)=\dfrac{t}{\bmbar h(t)}.
\end{equation*}
In particular,  the $A$-sequence and the function $h(t)$ determine each other. Furthermore, 
\begin{equation*}\label{TrivialAseq}
    A(t)=1 \iff h(t)=t
\end{equation*}
and the $Z$-sequence is determined by:
\begin{align*}\label{Zseq}
 d(t)&=\dfrac{d(0)}{1-tZ(h(t))}\\
 \iff d(t)-d(0)&=td(t)Z(h(t))\\
\iff
Z(t)&=\dfrac{d(\bmbar h(t))-d(0)}{\bmbar h(t)d(\bmbar h(t))} 
\end{align*}
\end{theorem}

The $A$- and $Z$-sequences for a  product of Riordan arrays have been completely determined as follows:

\begin{theorem}\label{thm:AZproducts}\cite[Theorem 3.3, Theorem 3.4]{HeSprugnoli2009} For the Riordan arrays $R_i=\mathcal {R}(d_i(t),h_i(t)), i=1,2$ with respective $A$-sequences and $Z$-se\-quences $A_i(t), Z_i(t), i=1,2,$ the $A$-sequence $A(t)$ and the $Z$-sequence $Z(t)$  of the product $R=R_1 R_2$ are given by 
\begin{equation*}\label{Aproduct} A(t)=A_2(t) A_1\left(\dfrac{t}{A_2(t)}\right) \end{equation*}
and 
\begin{equation*}\label{Zproduct} Z(t)=
\left(1-\dfrac{t}{A_2(t)}Z_2(t)\right)Z_1\left(\dfrac{t}{A_2(t)}\right)+  A_1\left(\dfrac{t}{A_2(t)}\right) Z_2(t). \end{equation*}
\end{theorem}

We will also need the following result, obtained by Lagrange Inversion \cite{StanEC2}.
\begin{theorem}\label{thm:InvRiordanElements}\cite[Theorem~2.2]{LuzonMerliniMoronSprugnoli2012LAA} 
Let $R=\R(d(t),h(t))$ be a Riordan array.  Then the inverse Riordan array $R^{-1}$ has $(n,k)$-entry equal to 
\[[t^{n-k}]\ \dfrac{h'(t)}{d(t) \left(\frac{h(t)}{t}\right)^{n+1}}.\]
\end{theorem}

For more on Riordan arrays and the Riordan group we refer the reader to \cite{shapiro1991riordan}, as well as the papers of Sprugnoli, e.g. \cite{SprugnoliFundThm1994}. 
\subsection{Examples}
\label{subsec:examples}
Obtaining combinatorial identities using the FTRA is straightforward. Here, we give three examples, with combinatorial proofs.  The first identity, Identity~\ref{id:fibo},  and its proof, are well known. 
To the best of our knowledge, the other two examples and Identities~\ref{id:2} and~\ref{id:riordanMotzkin} do not appear in the literature.

\begin{example}\label{ex:Fib-bisection}
Let $F_i$ denote the $i$th Fibonacci number. 
Let $(f_n)_{n\geq 0}$ be the sequence of squares of Fibonacci numbers, i.e., $$(f_n)_{n\geq 0} = (F_{n}^{2})_{n=0}^{\infty}=(1,1,4,9,25,64,169, \ldots ),$$ 
and let $(g_n)_{n\geq 0}$ denote the even-indexed Fibonacci numbers, i.e., $$(g_n)_{n\geq 0} = (F_{2n})_{n=0}^\infty =(1,2,5,13,34,89,233,\ldots).$$

The sequences $(f_n)_{n\geq 0}$ and $(g_n)_{n\geq 0}$ have generating functions 
$$f(t) = \frac{t(1-t)}{(1+t)(1-3t+t^2)} \quad \text{ and } \quad g(t) = \frac{t(1-t)}{(1-3t+t^2)},$$
respectively. Observing that $g(t) = (1+t)f(t)$, the FTRA gives that 
$D=\mathcal{R}(1+t,t)(f_n)_{n\geq 0} = (g_n)_{n\geq 0}$. 
In matrix form, this is
$$\begin{bmatrix}
1 & 0 & 0 & 0 & 0 & \cdots \\
1 & 1 & 0 & 0 & 0 & \cdots \\
0 & 1 & 1 & 0 & 0 & \cdots \\
0 & 0 & 1 & 1 & 0 & \cdots \\
0 & 0 & 0 & 1 & 1 & \cdots \\
\vdots & \vdots & \vdots & \vdots & \vdots & \ddots\\
\end{bmatrix}  
\cdot 
\begin{bmatrix}
1 \\
1 \\
4 \\
9 \\
25 \\
\vdots\\
\end{bmatrix}  
=
\begin{bmatrix}
1 \\
2 \\
5 \\
13 \\
34 \\
\vdots\\
\end{bmatrix}, 
$$
which gives the following well-known identity (see, for instance, \cite[Section 1.5]{benjamin}).
\begin{identity}
\label{id:fibo}
For $n\geq 0$, 
\begin{equation}
\label{eq:fibo}
F_{n}^{2}+F_{n+1}^{2}=F_{2n+1}.    
\end{equation}
\end{identity}
A combinatorial proof (mentioned in \cite{benjamin}) is given below, for completeness.  
\begin{proof}
We use the fact that the Fibonacci number $F_{k}$ is the number of ways to tile a board of length $k-1$ with squares and dominos \cite{StanEC1}. Thus, the right-hand side of \eqref{eq:fibo} is the number of ways to tile such a board of size $2n$. If the middle two cells with numbers $n$ and $n+1$ are tiled with a domino, we can tile the rest of the board in $F_{n}^{2}$ ways. If they are not, we have $F_{n+1}^{2}$ ways to tile the entire board.
\end{proof}
\end{example}

We record the following general proposition, obtained from the FTRA.

\begin{proposition}
Let $(f_n)_{n\geq 0}$ be a sequence with generating function $f(t)$ of the form $1/p(t)$ where $p(t)$ is a polynomial. 
Let $(g_n)_{n\geq 0}$ be a sequence with $g_0 \neq 0$ and generating function $g(t)$.
Then 
\begin{equation}
\label{eq:ex2}
\sum_{i=0}^n f_{n-i}\cdot \sum_{j=0}^m p_jg_{i-j} = g_n,
\end{equation}
with the understanding that $g_k=0$ whenever $k<0$.
\end{proposition}

\begin{proof}
Suppose we have a Riordan array $\mathcal{R}(d(t),t)$  which transforms $(f_n)_{n\ge 0}$ to $(g_n)_{n\geq 0}$, i.e., $\R(d(t),t))(f_i)_{i=0}^\infty = (g_i)_{i=0}^\infty$. 
Then FTRA tells us that 
\[ d(t) = g(t)/f(t)=g(t)p(t),\]
and the result follows.
\end{proof}

As a specific example, we relate the Fibonacci numbers to Catalan numbers. 
Recall \cite{StanCatalan2015} that the $i$th Catalan number is given by $C_i = \frac{1}{i+1}\binom{2i}{i}, i\ge 0$. 
The sequence of Catalan numbers is thus $(1,1,2,5,14,$ $42,132,429,\ldots )$.
\begin{example}\label{ex:Fib-to-shifted-Cat}

Let $(\bar{C}_n)_{n=0}^\infty$ denote the shifted Catalan sequence without the leading one. That is, the sequence $(\bar{C}_n)_{n=0}^\infty = (1,2,5,14,\ldots )$, with $n$th term given by $\bar{C}_n=C_{n+1}, n\ge 0$.  The generating function for $(\bar{C}_n)_{n= 0}^\infty$ is thus 
\[\bar{C}(t):=\frac{1-2t-\sqrt{1-4t}}{2t^2}.\]
(Note that the generating function for $(C_n)_{n\ge 0}$ 
is $t\bar{C}(t)+1=(1-\sqrt{1-4t})/(2t)$, see Section~\ref{subsec:Catalan}.
Since the Fibonacci sequence has generating function $F(t) = 1/p(t)$ with $p(t) = 1-t-t^2$, the Riordan array transforming $(F_i)_{i=0}^\infty$ to $(\bar{C}_i)_{i=0}^\infty$ is 

$$\R(\bar{C}(t)\cdot(1-t-t^2),t) = \R\left(\dfrac{\bar{C}(t)}{F(t)}, t\right) = \begin{bmatrix}
1 & 0 & 0 & 0 & 0 & 0 & \cdots \\
1 & 1 & 0 & 0 & 0 & 0 & \cdots \\
2 & 1 & 1 & 0 & 0 & 0 & \cdots \\
7 & 2 & 1 & 1 & 0 & 0 & \cdots \\
23 & 7 & 2 & 1 & 1 & 0 & \cdots \\
76 & 23 & 7 & 2 & 1 & 1 & \cdots \\
\vdots & \vdots & \vdots & \vdots & \vdots & \vdots & \ddots\\
\end{bmatrix},  
$$
where the $k$th column has $k$ zeros followed by the sequence $(\bar{C}_{i}-\bar{C}_{i-1}-\bar{C}_{i-2})_{i=0}^\infty$, with the understanding that $\bar{C}_i =0$ when $i<0$. This is the sequence  \cite[A067324]{OEIS} with two initial ones.

Equation~\ref{eq:ex2} simplifies to the following identity.

\begin{identity}
\label{id:2}
For $n\geq 3$, 
\begin{equation}
\label{eq:cat_fibo}
\sum\limits_{i=3}^{n} [\bar{C}_{i}-(\bar{C}_{i-1}+\bar{C}_{i-2})]F_{n+1-i} = \bar{C}_{n}-F_{n}-F_{n-1}.    
\end{equation}
\end{identity}

A combinatorial proof of Identity~\ref{id:2}  follows.

\begin{proof}
It is well known that $\bar{C}_{n}$ counts the number of Dyck paths, that is, paths from $(0,0)$ to $(2n,0)$, consisting of steps $U=(1,1)$ and $D=(1,-1)$, which never go below the x-axis. On the other hand, $F_{n+1}$ is the number of Dyck paths from $(0,0)$ to $(2n,0)$ consisting only of $UD$ and $UUDD$ segments. Thus, the right-hand side, $\bar{C}_{n}-F_{n}-F_{n-1}=\bar{C}_{n}-F_{n+1}$, is the number of the remaining Dyck paths which contain at least one segment different from these two. To interpret the left-hand side, note that the number of Dyck paths of length $2i$, which do not end in $UD$ or $UUDD$, is precisely $\bar{C}_{i}-(\bar{C}_{i-1}+\bar{C}_{i-2})$, since the number of paths ending in $UD$ is $\bar{C}_{i-1}$ and the number of paths ending in $UUDD$ is $\bar{C}_{i-2}$. Thus, the $i$th term of the sum on the left is the number of Dyck paths from $(0,0)$ to $(2n,0)$, such that they have a prefix ending at $(2i,0)$, not finishing with $UD$ or $UUDD$, but the rest of the path, from $(2i,0)$ to $(2n,0)$, is comprised of these two segments. It remains just to note that the possible values of $i$ are from $3$ to $n$.
\end{proof}
\end{example}

\begin{example}\label{ex:Riordan-to-Motzkin}
Consider the Motzkin numbers $(M_n)_{n=0}^\infty = (1,1,2,4,9,$ $21,51, \ldots)$ and the Riordan numbers $(R_n)_{n=0}^\infty = (1,0,1,1,3,$ $6,15,36,\ldots)$, along with their respective generating functions $M(t)$ and $R(t)$ (see Table~\ref{table:sequences}).

One can check that 
\begin{equation}
\label{eq:ex3d(t)}
\dfrac{M(t)}{R(t)}= tM(t)+1,
\end{equation}

and that the Riordan array $\R(tM(t)+1,t)$ below maps the sequence of Riordan numbers to the sequence of Motzkin numbers:
$$\begin{bmatrix}
1 & 0 & 0 & 0 & 0 & 0 & \cdots \\
1 & 1 & 0 & 0 & 0 & 0 & \cdots \\
1 & 1 & 1 & 0 & 0 & 0 & \cdots \\
2 & 1 & 1 & 1 & 0 & 0 & \cdots \\
4 & 2 & 1 & 1 & 1 & 0 & \cdots \\
9 & 4 & 2 & 1 & 1 & 1 & \cdots \\
\vdots & \vdots & \vdots & \vdots & \vdots & \vdots & \ddots\\
\end{bmatrix}  
\cdot 
\begin{bmatrix}
1 \\
0 \\
1 \\
1 \\
3 \\
6 \\
\vdots\\
\end{bmatrix}  
=
\begin{bmatrix}
1 \\
1 \\
2 \\
4 \\
9 \\
21\\
\vdots\\
\end{bmatrix}.
$$

Equation~(\ref{eq:ex3d(t)}) is equivalent to
\[tM(t)R(t)= M(t)-R(t), \]
which gives the following identity.
\begin{identity}
\label{id:riordanMotzkin}
For $n\geq 1$, 
\begin{equation}
\label{eq:riordan}
\sum\limits_{i=0}^{n-1} R_{i-1}M_{n-i} = M_{n}-R_{n}.    
\end{equation}
\end{identity}

A combinatorial proof  follows.

\begin{proof}
Recall that a Motzkin path of length $n$ is a lattice path in the plane from $(0, 0)$ to $(n, 0)$ never going below the $x$-axis, which consists of up steps
$U = (1, 1)$, down steps $D = (1,-1)$, and horizontal steps $H = (1, 0)$. The Motzkin number $M_{n}$ counts the number of Motzkin paths of length $n$. On the other hand, the Riordan number $R_{n}$ counts the Motzkin paths of length $n$ with no horizontal steps of height $0$ \cite[Section 2]{chen}.
Therefore, the right-hand side of \eqref{eq:riordan} is the number of Motzkin paths of length $n$, which have at least one horizontal step of height $0$. Let the first such step be step $i$. Before this step, we must have a path of length $i-1$ without horizontal steps of height $0$. The number of these paths is $R_{i-1}$. After this step, we can have any Motzkin path of length $n-i$ and thus we have $M_{n-i}$ such paths. Summing over the possible values of $i$ gives the left-hand side.  
\end{proof}

\end{example}

\section{Sums of Riordan arrays}\label{sec:Sums}
Although Riordan arrays form a group under multiplication, the sum of Riordan arrays is not necessarily a Riordan array, as has already been observed in the literature, e.g. \cite[p. 172]{MERLINI2017160}. 
In this section we examine the summations of Riordan arrays more closely. 
We first determine when the sum of Riordan arrays yields a Riordan array (Theorem~\ref{thm:R+S}). 
We then consider the case of summing arbirary Riordan arrays and show that it satisfies a recurrence (Theorem~\ref{thm:Eric}).

\subsection{When addition of Riordan arrays is closed}

The following theorem specifies exactly when the sum of Riordan arrays is again a Riordan array.

\begin{theorem}
\label{thm:R+S}
Let $R = \mathcal{R}(d_R(t), h_R(t))$ and $S = \mathcal{R}(d_S(t), h_S(t))$ be two Riordan arrays. Then $R+S$ is a Riordan array if and only if $h_R(t)=h_S(t)=h(t)$, say, and $d_R(0) + d_S(0) \ne 0$. In that case $R+S$ is the Riordan array $\mathcal{R}(d_R(t) + d_S(t), h(t))$.
\end{theorem}



\begin{proof}  Suppose $h_R(t)=h_S(t)=h(t)$ and $d_R(0) + d_S(0) \ne 0$. By definition, 
\[R_{n,c} = [t^n]\{ d_R(t) h(t)^c \}\] and 
\[S_{n,c} = [t^n]\{ d_S(t) h(t)^c \}.\]
Hence
\begin{align*}
R_{n,c} + S_{n,c} &= [t^n]\{ d_R(t) h(t)^c + d_S(t)h(t)^c\} \\
&= [t^n]\{ (d_R(t) + d_S(t)) h(t)^c \}.
\end{align*}
Since $d_R(0) + d_S(0) \ne 0$, this gives that $R + S$ is a Riordan array.

Conversely, suppose $R+S$ is a Riordan array. Hence there are series $d_{R+S}(t), h_{R+S}(t)$ such that $R+S=\mathcal{R}(d_{R+S}(t), h_{R+S}(t))$.  From Equation~\eqref{d} it is clear that $d_{R+S}(t)=d_R(t)+d_S(t).$  By Theorem~\ref{thm:FTRA},  for any power series $f(t)$, 
\begin{enumerate}
\item
$R$ transforms $f(t)$ to $d_R(t) f(h_R(t))$, 
\item
$S$ transforms $f(t)$ to $d_S(t) f(h_S(t))$,
\item
$R+S$ transforms $f(t)$ to $d_{R+S}(t) f(h_{R+S}(t))$.
\end{enumerate}
Thus
\begin{equation*} d_{R+S}(t) f(h_{R+S}(t))=d_R(t) f(h_R(t))+d_S(t) f(h_S(t)).
\end{equation*}
Equivalently, since   Equation~\eqref{d} implies that $d_{R+S}(t)=d_R(t)+d_S(t)$, we have 
\begin{equation*}  f(h_{R+S}(t))=\dfrac{d_R(t) f(h_R(t))+d_S(t) f(h_S(t))}{d_R(t)+d_S(t)}.
\end{equation*}


Applying this to the polynomial $f(t)=t^k$, $k= 1,2,$ and suppressing the $t$'s for clarity, gives\footnote{The authors thank Matt Hudelson for completing this proof.}
\begin{equation*}  h_{R+S}=\dfrac{d_R h_R+d_S h_S}{d_R+d_S},\quad h_{R+S}^2=\dfrac{d_R h_R^2+d_S h_S^2}{d_R+d_S}.
\end{equation*}
This in turn implies 
\begin{align*} &(d_R h_R^2+d_S h_S^2 )(d_R+d_S)=(d_R h_R+d_S h_S)^2\\
&\iff d_Rd_S(h_R^2+h_S^2)= 2  d_Rd_S h_R h_S\\
&\iff (h_R-h_S)^2=0,
\end{align*}
since $d_Rd_S\ne 0,$ thereby finishing the proof.
\end{proof}


\begin{remark}
Let $G_{h(t)}$ denote a subgroup of Riordan group consisting of the Riordan arrays $\R(d(t),h(t))$ for a fixed $h(t)$.  
Although the array of zeroes $\R(0,h(t))$ is not a Riordan array, the set   $G_{h(t)}\cup \{\R(0,h(t)) \}$ admits a ring structure under matrix addition and multiplication.
\end{remark}

We can now give an analogue of Theorem~\ref{thm:AZproducts}, characterising the $A$- and $Z$-sequences formed by the sum of Riordan arrays with a common $h(t)$.

\begin{theorem}\label{ZSeqR+S}
Let $R = \mathcal{R}(d_R(t), h(t))$ and $S = \mathcal{R}(d_S(t), h(t))$ be two Riordan arrays with the same $h(t)$, such that $d_R(0) + d_S(0) \ne 0$. Then $R$, $S$, and $R+S$ have the same A-sequence $A(t)=\dfrac{t}{\bmbar h (t)}.$ The Z-sequence of $R+S$ is given by 
\begin{equation}\label{eqn:ZofR+Sv2}
Z_{R+S}(t)=\dfrac{d_R(0)Z_R(t)+d_S(0)Z_S(t) -d_{R+S}(0) \bmbar h(t) Z_R(t) Z_S(t) }{d_{R+S}(0)  -\bmbar h(t) d_R(0)Z_S(t) -\bmbar h(t) d_S(0)Z_R(t) };
\end{equation}

This expression can also be written in several equivalent ways as follows. 
We have 
\begin{equation}\label{ZofR+Sv1}
Z_{R+S}(t)=\dfrac{[d_R(0)Z_R(t)+d_S(0)Z_S(t)] A(t) -td_{R+S}(0)  Z_R(t) Z_S(t) }{d_{R+S}(0)A(t) -t [d_R(0)Z_S(t)+d_S(0)Z_R(t) ] };\end{equation}
\begin{equation}\label{eqn:Sheila2}Z_{R+S}(h(t))=\frac{d_R(t)+d_S(t)-d_R(0)-d_S(0)}{t(d_R(t)+d_S(t))};\end{equation}
and finally a symmetric expression as a weighted average of the individual $Z$-sequences:
\begin{equation}\label{eqn:Sheila3}Z_{R+S}(t)=\frac{d_R(\bmbar h (t))Z_R(t)+d_S(\bmbar h (t))Z_S(t)}{d_R(\bmbar h (t))+d_S(\bmbar h (t))}
\end{equation}
\end{theorem}
\begin{proof} By Theorem \ref{thm:He-Sprugnoli} and Theorem \ref{thm:R+S}, since $R$, $S$, and $R+S$ have the same $h(t)$, it follows that $R$, $S$, and $R+S$ have the same A-sequence. By Theorem \ref{thm:He-Sprugnoli}, this is $A(t)=\dfrac{t}{\bmbar h (t)}$.


From Theorem \ref{thm:He-Sprugnoli} we also know that 
\[ d_R(\bmbar h(t))=\dfrac{d_R(0)}{1-\bmbar h(t) Z_R(t)},\   d_S(\bmbar h(t))=\dfrac{d_S(0)}{1-\bmbar h(t) Z_S(t)}, \]
and 
\begin{equation}\label{dRplusS-of-hbar}d_{R+S}(\bmbar h(t))=\dfrac{d_{R+S}(0)}{1-\bmbar h(t) Z_{R+S}(t)}.\end{equation}

Since $d_{R+S}=d_R+d_S$ by Theorem \ref{thm:R+S}, we have (suppressing the $t$'s for clarity) that 

\begin{align*}
\dfrac{d_{R+S}(0)}{1-\bmbar h Z_{R+S}}
& = \dfrac{d_R(0)}{1-\bmbar h Z_R}+\dfrac{d_S(0)}{1-\bmbar h Z_S}
 \\
& = \dfrac{d_{R+S}(0) -\bmbar h \ d_R(0)Z_S -\bmbar h\ d_S(0)Z_R }{(1-\bmbar h Z_R ) (1-\bmbar h Z_S)}.
\end{align*}
Hence we have 
$$
    d_{R+S}(0)(1-\bmbar h  Z_R)(1-\bmbar h Z_S)
    =\left[ d_{R+S}(0) -\bmbar h (d_R(0)Z_S+d_S(0)Z_R )\right] (1-\bmbar h Z_{R+S}).
$$
The left-hand side equals 
$$ d_{R+S}(0) -\bmbar h \left[d_{R+S}(0) \left(Z_R+Z_S\right) -d_{R+S}(0)\bmbar h Z_R Z_S\right].$$
The right-hand side equals 
$$
    d_{R+S}(0) -\bmbar h (d_R(0)Z_S+d_S(0)Z_R ) + [d_{R+S}(0) -\bmbar h (d_R(0)Z_S+d_S(0)Z_R )] (-\bmbar h Z_{R+S}).
$$
It follows that 
$$ -\bmbar h \left[d_R(0)Z_R+d_S(0)Z_S -d_{R+S}(0) \bmbar h Z_R Z_S \right]
=[d_{R+S}(0) -\bmbar h (d_R(0)Z_S+d_S(0)Z_R ) ] (-\bmbar h Z_{R+S})
$$
again using $d_{R+S}=d_R+d_S,$ and hence 
$$
    d_R(0)Z_R+d_S(0)Z_S -d_{R+S}(0) \bmbar h Z_R Z_S  
    =\left[d_{R+S}(0) -\bmbar h \left(d_R(0)Z_S+d_S(0)Z_R \right) \right]  Z_{R+S}.  
$$
This gives (\ref{eqn:ZofR+Sv2}). 
Now using $\bar{h}(t)=t/A(t)$ gives \eqref{ZofR+Sv1}.

Replacing $t$ with $h(t)$ in ~\eqref{dRplusS-of-hbar} and rearranging gives the expression \eqref{eqn:Sheila2} for the composition $Z_{R+S} \circ h$ in terms of only the $d(t)$'s.
Finally, since $d(t)-d(0)=td(t)Z(h(t))$ by Theorem \ref{thm:He-Sprugnoli}, we can also write the composition $Z_{R+S} \circ h$ as a symmetric expression in $d(t)$ and $Z(t)$:
\begin{equation*}\label{Rachel}Z_{R+S}(h(t))=\frac{d_R(t)Z_R(h(t))+d_S(t)Z_S(h(t))}{d_R(t)+d_S(t)},
\end{equation*}
from which Equation~\eqref{eqn:Sheila3} follows.
\end{proof}

An important special case is the following:
\begin{corollary}\label{cor:Z_for_h(t)_is_t}  Let $h(t)=t=\bmbar h (t),$ and thus $A(t)=1.$
 
Then equations~\eqref{eqn:ZofR+Sv2}, \eqref{eqn:Sheila2} and \eqref{eqn:Sheila3} become respectively
\begin{align}\label{ZofR+Sv3}
Z_{R+S}(t)&=\dfrac{d_R(0)Z_R(t)+d_S(0)Z_S(t) -d_{R+S}(0) t\, Z_R(t) Z_S(t) }{d_{R+S}(0)  -t \ d_R(0)Z_S(t) -t\ d_S(0)Z_R(t) }\\
&=\frac{d_R(t)+d_S(t)-d_R(0)-d_S(0)}{t(d_R(t)+d_S(t))}\\
&=\frac{d_R(t)Z_R(t)+d_S(t)Z_S(t)}{d_R(t)+d_S(t)}.
\end{align}
\end{corollary}

\begin{example}
\label{ex:addRiordanAZsequence}
Consider the Riordan arrays $\Pas = \mathcal{R}\left(\frac{1}{1-t}, \frac{t}{1-t}\right)$, and
$S = \mathcal{R}\left(\frac{1}{1-2t}, \frac{t}{1-t}\right)$.
Their $A$- and $Z$-sequences are $A_{\Pas}=(1,1,0,0,0,\dots)$, $Z_{\Pas}=(1,0,0,0,\dots)$, $A_S=(1,1,0,0,0,\dots)$, and $Z_S=(2,0,0,0,\dots)$, or, equivalently, $A_{\Pas}(t)=A_S(t)=1+t$, $Z_{\Pas}(t)=1$, and $Z_S(t)=2$.

We have 
$$\Pas=\begin{bmatrix}
1 & 0 & 0 & 0 & 0 & \cdots \\
1 & 1 & 0 & 0 & 0 & \cdots \\
1 & 2 & 1 & 0 & 0 & \cdots \\
1 & 3 & 3 & 1 & 0 & \cdots \\
1 & 4 & 6 & 4 & 1 & \cdots \\
\vdots & \vdots & \vdots & \vdots & \vdots & \ddots\\
\end{bmatrix}, \hspace{5 mm} S=\begin{bmatrix}
1 & 0 & 0 & 0 & 0 & \cdots \\
2 & 1 & 0 & 0 & 0 & \cdots \\
4 & 3 & 1 & 0 & 0 & \cdots \\
8 & 7 & 4 & 1 & 0 & \cdots \\
16 & 15 & 11 & 5 & 1 & \cdots \\
\vdots & \vdots & \vdots & \vdots & \vdots & \ddots\\
\end{bmatrix}.
$$

Their sum is 
$$\Pas+S=\mathcal{R}\left(\frac{2-3t}{1-3t+2t^2}, \frac{t}{1-t}\right)=\begin{bmatrix}
2 & 0 & 0 & 0 & 0 & \cdots \\
3 & 2 & 0 & 0 & 0 & \cdots \\
5 & 5 & 2 & 0 & 0 & \cdots \\
9 & 10 & 7 & 2 & 0 & \cdots \\
17 & 19 & 17 & 9 & 2 & \cdots \\
\vdots & \vdots & \vdots & \vdots & \vdots & \ddots\\
\end{bmatrix}
$$

Using Theorem \ref{ZSeqR+S}, we determine that the $A$- and $Z$-sequences for the Riordan array $\Pas+S$ will be $A_{\Pas+S}(t)=1+t$, or $A_{\Pas+S}=(1,1,0,0,0,\dots)$, and $$Z_{\Pas+S}(t)=\dfrac{[1\cdot 1+1\cdot 2] (1+t) -t \cdot 2 \cdot 1 \cdot 2 }{2(1+t) -t [1\cdot 2+1\cdot 1 ] }=\frac{3-t}{2-t},$$
or $Z_{\Pas+S}=\left(\frac{3}{2},\frac{1}{4},\frac{1}{8},\frac{1}{16},\dots\right)$.
\end{example}

In the next example, we demonstrate how the summation of Riordan arrays can be used to obtain combinatorial identities. 

\begin{example}
Let $R = \R\left(\frac{-1}{1-t^2}, \frac{t}{1-t^2}\right)$ and 
$S = \R\left(\frac{2}{(t^2-1)^2}, \frac{t}{1-t^2}\right)$.
We have that 
$$R =\R\left(\frac{-1}{1-t^2}, \frac{t}{1-t^2}  \right) = 
\begin{bmatrix}
-1 & 0 & 0 & 0 & 0 &  \cdots\\
0 & -1 & 0 & 0 & 0 &  \cdots\\
-1 & 0 & -1 & 0 & 0 &  \cdots\\
0 & -2 & 0 & -1 & 0 &  \cdots\\
-1 & 0 & -3 & 0 & -1 & \cdots \\
\vdots & \vdots & \vdots & \vdots & \vdots & \ddots
\end{bmatrix},$$

with general term given by
$$R_{n,k} = -\frac{1+(-1)^{n+k}}{2}\binom{(n+k)/2}{k}.$$

We also have that 
$$S =\R\left(\frac{2}{(1-t^2)^2}, \frac{t}{1-t^2}  \right) = 
2\cdot \begin{bmatrix}
1 & 0 & 0 & 0 & 0 &  \cdots\\
0 & 1 & 0 & 0 & 0 &  \cdots\\
2 & 0 & 1 & 0 & 0 &  \cdots\\
0 & 3 & 0 & 1 & 0 &  \cdots\\
3 & 0 & 4 & 0 & 1 & \cdots \\
\vdots & \vdots & \vdots & \vdots & \vdots & \ddots
\end{bmatrix},$$

with general term given by 
$$ 
S_{n,k} = (1 + (-1)^{n+k})\binom{(n+k+2)/2}{k+1}.
$$

The sum $R+S$ is then computed to be
$$R + S =\R\left(\frac{t^2+1}{(t^2-1)^2}, \frac{t}{1-t^2}  \right) = 
\begin{bmatrix}
1 & 0 & 0 & 0 & 0 &  \cdots\\
0 & 1 & 0 & 0 & 0 &  \cdots\\
3 & 0 & 1 & 0 & 0 &  \cdots\\
0 & 4 & 0 & 1 & 0 &  \cdots\\
5 & 0 & 5 & 0 & 1 & \cdots \\
\vdots & \vdots & \vdots & \vdots & \vdots & \ddots
\end{bmatrix},$$

with general term given by 
$$ (R+S)_{n,k} = \frac{1+(-1)^{n-k}}{2} \cdot \frac{n+1}{k+1}\binom{(n+k)/2}{k}
$$

The row sums of $R$ are the negative Fibonacci numbers $(-F_i)_{i=0}^\infty = (-1,-1$, $-2$, $-3$, $-5,...)$, while the row sums of $S$ are given by $2 (a_i)_{i=0}^\infty = 2F_{i+1} - 1-(-1)^{i+1}$, where $(a_i)$ is the sequence \cite[A052952]{OEIS}.
The row sums of $R+S$ give the associated Mersenne numbers $1,1,4,5,11,16,29, \ldots$ (\cite[A001350]{OEIS}), with the $i$th term given by $L_{i+1}- 1 - (-1)^{i+1}$.

As a consequence, we have
$$ -F_i + 2F_{i+1} - 1 - (-1)^{i+1} = L_{i+1} -1 - (-1)^{i+1},$$
from which we obtain the following identity  \cite[A000032]{OEIS}.
\begin{identity}
\label{id:fibLuc}
$$ 2F_{i-1} + F_i = L_{i+1} $$
\end{identity}

We observe that $R\cdot (1,2,3,4,5, \ldots )^t  = (-1, -2, -4, -8, -15, -28, -51, \ldots )^t$. 
This is the negative of a shift of the sequence \cite[A029907]{OEIS}, and its $n$th term is given by $-\frac{1}{5}((n+5)F_{n} + 2(n+1)F_{n-1})$.
Using the FTRA, its generating function can be verified to be $d_R(t)f(h_R(t)) = (1-t^2)/(1-t-t^2)^2$ with $f(t) = 1/(1-t)^2$. 

Similarly, we observe that  $S\cdot (1,2,3,4,5, \ldots )^t  = 2(1, 2, 5, 10, 20, 38, 71,\ldots )^t$, which is twice the convolution of the Fibonacci sequence with itself, i.e., $2 \sum_{k=0}^{n}F_k F_{n-k}$.
Furthermore, it is easy to check that 
$(R+S)$ transforms the sequence $(1,2,3,4,5, \ldots)$ to the sequence $((n+1)F_{n})_{n= 0}^\infty$,
which is \cite[A045925]{OEIS}.  

Thus we have that 
$$-\frac{1}{5}((n+5)F_{n} + 2(n+1)F_{n-1}) + 2 \sum_{k=0}^{n}F_k F_{n-k} = (n+1)F_{n},$$
which simplifies to the following expression for the convolution of Fibonacci numbers with themselves.
\begin{identity}
\label{id:fibConv}
For $n\geq 0$,
$$\sum_{k=0}^{n}F_kF_{n-k} = \frac{(3n+5)F_{n} + (n+1)F_{n-1}}{5}.$$
\end{identity}
This identity is closely related to an identity given in \cite[A001629]{OEIS}.
\end{example}

\subsection{Sums of arbitrary Riordan arrays}\label{sec:sum-arb}

Let $R = \mathcal{R}(d_R(t), h_R(t))$ and $S = \mathcal{R}(d_S(t), h_S(t))$ be two Riordan arrays where $h_R(t)$ may not equal  $h_S(t)$. The matrix sum $R + S$ is not necessarily a Riordan array but does have some interesting properties in its own right. We call such an array a \textit{Riordan sumray} and use the notation $$R+S=\mathcal{R}(d_R(t), d_S(t), h_R(t), h_S(t))$$ to specify it. Note that the $(n,k)$-element of $\mathcal{R}(d_R(t), d_S(t), h_R(t), h_S(t))$ is the coefficient on $t^n$ in $d_R(t)(h_R(t))^k + d_S(t)(h_S(t))^k$. 

\begin{example}
\label{ex:sumray}
From the Riordan arrays $R = \mathcal{R}\left(\frac{1}{1-t}, \frac{t}{1-t}\right)$ and $S = \mathcal{R}\left(\frac{1}{1-t}, \frac{2t}{1-t}\right)$ we have the Riordan sumray:
$$
\mathcal{R}\left(\frac{1}{1-t}, \frac{1}{1-t}, \frac{t}{1-t}, \frac{2t}{1-t}\right)=\begin{bmatrix}
2 & 0 & 0 & 0 & 0 & 0 & \cdots\\
2 & 3 & 0 & 0 & 0 & 0 & \cdots\\
2 & 6 & 5 & 0 & 0 & 0 & \cdots\\
2 & 9 & 15 & 9 & 0 & 0 & \cdots\\
2 & 12 & 30 & 36 & 17 & 0 & \cdots\\
2 & 15 & 50 & 90 & 85 & 33 & \cdots\\
\vdots & \vdots & \vdots & \vdots & \vdots & \vdots & \ddots
\end{bmatrix}
$$
\end{example}

Our goal in this subsection is to derive an interesting second order recurrence for an arbitrary sum $R+S$. 
We begin with some notation and preliminary lemmas.  
\begin{definition}\label{defn:DownshiftA} Let $A:=(a_n)_{n\ge 0}$ be any sequence.  We define  
$\tilde{A}$ to be the infinite matrix whose $(i,j)$-entry is 
$$\begin{cases} a_{i-j}, & i\ge j,\\ 0, & \text{otherwise}. \end{cases}$$ 
In particular, $\tilde{A}$ is a Toeplitz array; it is lower triangular with constant entries along the main diagonal and its parallel subdiagonals, where the $n$th term of the sequence fills the subdiagonal formed by $\{(i,j):i-j=n\}$. 
Note that the $j$th column of $\tilde{A}$ is simply the sequence $(a_n)_{n\ge 0}$ shifted down by $j$ spaces.
\end{definition}

\begin{remark}
Note that if the sequence from Definition \ref{defn:DownshiftA} is an $A$-sequence $A:=(a_0 \neq 0, a_1 ,a_2, ...)$, one has $\tilde{A}=\mathcal{R}(A(t),t).$  
\end{remark}

\begin{lemma}\label{lem:CommutingAS} Let $A:=(a_n)_{n\ge 0}$ and $B:=(b_n)_{n\ge 0}$ be any two sequences.  Then the matrices $\tilde{A}$ and $\tilde{B}$ commute.
\end{lemma}
\begin{proof} The $(i,j)$-entry of the product $\tilde{A}\tilde{B}$ is the convolution 
\[\sum_{k\ge 0} a_{i-k}b_{k-j}.\]
This is nonzero if and only if $i\ge k\ge j,$ and then it equals 
\[\sum_{k=j}^i a_{i-k}b_{k-j}=a_{i-j} b_0+a_{i-j-1} b_1 +\cdots + a_0b_{i-j}=\sum_{\stackrel{p,q\ge 0}{p+q=i-j}} a_p b_q.\]
Clearly the latter convolution is symmetric in $A$ and $B$.
\end{proof}

\begin{lemma}\label{lem:TruncatedRiordan} Let $R = \mathcal{R}(d_R(t), h_R(t))$  be a Riordan array, and let $A_R$ be its $A$-sequence.  Let $R^\sigma$ denote the array obtained from $R$ by deleting the first row and the first column.  Then
\begin{enumerate}
\item 
 $R^\sigma$ is the Riordan array 
$\mathcal{R}(d_R(t) h_R(t)/t, h_R(t))$, and hence has the same $A$-sequence $A_R$;
\item the matrix equation $R^\sigma=R\, \tilde{A}_R$ holds.
\end{enumerate}
\end{lemma}

\begin{proof} Note that the coefficient of $t^n$ in $(d_R(t) h_R(t)/t) \cdot 
h_R(t)^k$ equals the coefficient of $t^{n+1}$ in $d_R(t) h_R(t)^{k+1}$, which is precisely the 
$(n+1, k+1)$-entry of the Riordan array  $R = \mathcal{R}(d_R(t), h_R(t))$, from Equation~\eqref{d}.  But this is also the $(n,k)$-entry of $R^\sigma$ for $n,k\ge 0,$ and the first part of the first statement follows.  
The $A$-sequence is preserved since it depends only on $h_R(t)$ by Theorem~\ref{thm:He-Sprugnoli}. 

Let $d_{n,k}$ denote the $(n,k)$-entry of the Riordan array $R$. For the second statement, the definition of the $A$-sequence gives (Equation~\eqref{A}) for $n,k\ge 0$, 
\[d_{n+1,k+1}=\sum_{i\ge 0} d_{n,k+i} a_i=\sum_{ j\ge 0} d_{n,j} a_{j-k}.\]
But the right-hand side is the $(n,k)$-entry of the matrix product $R\, \tilde{A}_R$, by Definition~\ref{defn:DownshiftA}, while $d_{n+1,k+1}$ is the $(n,k)$-entry of $R^\sigma$.  The claim follows. 
\end{proof}

\begin{remark}
Part (1) of Lemma \ref{lem:TruncatedRiordan} also follows from Theorem \ref{thm:Shapiro} and Theorem \ref{thm:He-Sprugnoli} since
\begin{align*}
\mathcal{R}\left(d_R(t), h_R(t)\right)\mathcal{R}\left(A(t),t\right)&=
\mathcal{R}\left(d_R(t), h_R(t)\right)\mathcal{R}\left(\frac{t}{\bmbar{h_R}(t)},t\right)\\
&=\mathcal{R}\left(d_R(t)\frac{h_R(t)}{\bmbar{h_R}(h_R(t))}, h_R(t)\right)\\
&=\mathcal{R}\left(d_R(t)\frac{h_R(t)}{t}, h_R(t)\right)
\end{align*}
See also \cite[Proof of Theorem~5.3.1]{Sprugnoli2006} 
\end{remark}

We can now establish a recurrence relating the sum of any two arbitrary Riordan arrays.  This recurrence for Riordan sumrays is an analog of the typical $A$-sequence.

\begin{theorem}\label{thm:Eric} Let $R = \mathcal{R}(d_R(t), h_R(t))$ and $S = \mathcal{R}(d_S(t), h_S(t))$ be two Riordan arrays with respective $A$-sequences $A_R, A_S.$ Then the Riordan sumray $R + S$ satisfies the following recurrence.
\[ (R+S)^{\sigma\sigma}= (R + S)^\sigma (\tilde{A}_R + \tilde{A}_S) - (R + S) \tilde{A}_R \tilde{A}_S.\]
Equivalently, the entries $d_{n,k}$ of the sum $R+S$ satisfy the following recurrence:
     \[ d_{n+2,k+2} = \sum_{j = 0}^\infty d_{n+1,j+1} B_{j,k}+ \sum_{j=0}^\infty d_{n,k+j}C_{j,k}, \]
    where $B=\tilde{A}_R+\tilde{A}_S$, and $C=-\tilde{A}_R\tilde{A}_S$, and thus 
$B_{j,k}= (A_R)_{j-k}+(A_S)_{j-k}$  and $C_{j,k}=-\sum_ {\stackrel{p,q\ge 0}{p+q=j-k}} (A_R)_p (A_S)_q $.
\end{theorem}
\begin{proof} Note that, by using Part (2) followed by Part (1) of Lemma~\ref{lem:TruncatedRiordan},
$$R^{\sigma \sigma} = R^\sigma \tilde{A}_{R^\sigma} = R^\sigma\tilde{A}_{R}.$$ 
Similarly,
$$S^{\sigma \sigma} = S^\sigma \tilde{A}_{S^\sigma} = S^\sigma\tilde{A}_{S}.$$
Since  $(R+S)^\sigma = R^\sigma + S^\sigma$, we have:
    \begin{align*}
        (R + S)^{\sigma \sigma} &= R^{\sigma \sigma} + S^{\sigma \sigma} \\
        &= R^\sigma \tilde{A}_R + S^\sigma \tilde{A}_S \\
        &= (R^\sigma + S^\sigma)(\tilde{A}_R + \tilde{A}_S) - R^\sigma  \tilde{A}_S - S^\sigma  \tilde{A}_R \\
        &= (R^\sigma + S^\sigma)(\tilde{A}_R + \tilde{A}_S) - R \tilde{A}_R \tilde{A}_S - S \tilde{A}_S \tilde{A}_R \\
        &= (R + S)^\sigma (\tilde{A}_R + \tilde{A}_S) - (R + S) \tilde{A}_R \tilde{A}_S,
    \end{align*}
since $\tilde{A}_R$ and $\tilde{A}_S$ commute by Lemma~\ref{lem:CommutingAS}. 
    
Let $d_{n,k}, n,k\ge 0,$ denote the $(n,k)$-entry of the sum $R+S$.  Then the $(n,k)$-entries of $(R+S)^\sigma$ and $(R+S)^{\sigma\sigma}$ are respectively  
\[d_{n+1,k+1} \ \text{ and } \ d_{n+2,k+2} \]
for $n,k \geq 0$.
 Writing out the corresponding equation for the matrix entries gives, for $n,k\ge 0$: 
    \[ d_{n+2,k+2} = \sum_{j = 0}^\infty   d_{n+1,j+1} B_{j,k}+ \sum_{j=0}^\infty   d_{n,k+j}C_{j,k}, \]
    where $B=\tilde{A}_R+\tilde{A}_S$, and $C=-\tilde{A}_R\tilde{A}_S$.
%
\end{proof}

\begin{example}
Consider Example \ref{ex:sumray} within the context of this theorem. The A-sequences for the Riordan arrays which add to produce this Riordan sumray are $(1,1,0,0,\ldots)$ and $(2,1,0,0,\ldots)$. Therefore, the $B$ and $C$ matrices are:
$$
B=\begin{bmatrix}
3 & 0 & 0 & 0 & 0 &  \cdots\\
2 & 3 & 0 & 0 & 0 &  \cdots\\
0 & 2 & 3 & 0 & 0 &  \cdots\\
0 & 0 & 2 & 3 & 0 &  \cdots\\
0 & 0 & 0 & 0 & 2 &  \cdots\\
\vdots & \vdots &  \vdots & \vdots & \vdots & \ddots
\end{bmatrix}
\qquad 
C=\begin{bmatrix}
-2 & 0 & 0 & 0 & 0 &  \cdots\\
-3 & -2 & 0 & 0 & 0 & \cdots\\
-1 & -3 & -2 & 0 & 0 &  \cdots\\
0 & -1 & -3 & -2 & 0 &  \cdots\\
0 & 0 & -1 & -3 & -2 &  \cdots\\
\vdots & \vdots & \vdots & \vdots & \vdots & \ddots
\end{bmatrix}
$$

The Riordan sumray $R = \mathcal{R}\left(\frac{1}{1-t}, \frac{1}{1-t}, \frac{t}{1-t}, \frac{2t}{1-t}\right)$ then satisfies the second-order recurrence:
\[ R^{\sigma\sigma} = R^\sigma B + R C \]
\end{example}

The previous example illustrates the recurrence in terms of array multiplication. We can also use this recurrence to calculate specific entries $d_{n,k}$ as shown in the next example.

\begin{example}
\label{ex:addRiordanSumrays}
    Consider the Pascal array $\Pas=\mathcal{R}\left(\frac{1}{1-t},\frac{t}{1-t}\right)$ and the Shapiro array $\Sha=\mathcal{R}\left(\frac{1-\sqrt{1-4t}}{2t},\frac{1-2t-\sqrt{1-4t}}{2t}\right)$.
    $$\Pas=\begin{bmatrix}
    1 & 0 & 0 & 0 & 0 & \cdots \\
    1 & 1 & 0 & 0 & 0 & \cdots \\
    1 & 2 & 1 & 0 & 0 & \cdots \\
    1 & 3 & 3 & 1 & 0 & \cdots \\
    1 & 4 & 6 & 4 & 1 & \cdots \\
    \vdots & \vdots & \vdots & \vdots & \vdots & \ddots\\
    \end{bmatrix}, \hspace{5 mm} \Sha=\begin{bmatrix}
    1 & 0 & 0 & 0 & 0 & \cdots \\
    1 & 1 & 0 & 0 & 0 & \cdots \\
    2 & 3 & 1 & 0 & 0 & \cdots \\
    5 & 9 & 5 & 1 & 0 & \cdots \\
    14 & 28 & 20 & 7 & 1 & \cdots \\
    \vdots & \vdots & \vdots & \vdots & \vdots & \ddots\\
    \end{bmatrix},
    $$
    and 
    $$\Pas+\Sha=\begin{bmatrix}
    2 & 0 & 0 & 0 & 0 & \cdots \\
    2 & 2 & 0 & 0 & 0 & \cdots \\
    3 & 5 & 2 & 0 & 0 & \cdots \\
    6 & 11 & 8 & 2 & 0 & \cdots \\
    15 & 32 & 26 & 11 & 2 & \cdots \\
    \vdots & \vdots & \vdots & \vdots & \vdots & \ddots.\\
    \end{bmatrix}
    $$
    
    The $A$-sequences for $\Pas$ and $\Sha$ are $A_{\Pas}=(1,1,0,0,\dots)$ and $A_{\Sha}=(1,2,1,0,0,\dots)$. We check that the entries $d_{n+2,k+2}$ for $n,k \geq 0$ of $\Pas+\Sha$ satisfy the formula given by Theorem \ref{thm:Eric}.  For instance, the following calculation gives us the $d_{3,2}$ entry:
    \begin{align*}
        8=d_{3,2} &= \sum_{j=0}^{\infty}d_{2,j+1}B_{j,0} + \sum_{j=0}^{\infty} d_{1,j}C_{j,0}\\
        &= d_{2,1}B_{0,0}+d_{2,2}B_{1,0}+d_{1,0}C_{0,0}+d_{1,1}C_{1,0}\\
        &= 5(2) + 2(3) + 2(-1) + 2(-3),
    \end{align*}
    where
    \begin{align*}
        &B_{0,0} = (\tilde{A}_{\Pas})_{0,0}+(\tilde{A}_{\Sha})_{0,0} & & C_{0,0} = -\left((\tilde{A}_{\Pas})_{0,0}(\tilde{A}_{\Sha})_{0,0}\right)\\
        &B_{1,0} = (\tilde{A}_{\Pas})_{1,0} + (\tilde{A}_{\Sha})_{1,0} & & C_{1,0} = -\left((\tilde{A}_{\Pas})_{1,0}(\tilde{A}_{\Sha})_{0,0} + (\tilde{A}_{\Pas})_{1,1}(\tilde{A}_{\Sha})_{1,0}\right).
    \end{align*}
\end{example}

When $h_R(t) = h_S(t)$, $R + S$ is itself a Riordan array with $\tilde{A}_R =\tilde{A}_S=\tilde{A}$, and the recurrence from Theorem \ref{thm:Eric} specializes to 
\begin{align*}
    (R+S)^{\sigma\sigma} &= 2(R + S)^\sigma (\tilde{A}) - (R + S) \tilde{A}^2\\
    &= 2(R + S) (\tilde{A})(\tilde{A}) - (R + S) \tilde{A}^2\\
    &=(R+S) \tilde{A}^2,
\end{align*}
as expected. Revisiting Example \ref{ex:addRiordanAZsequence}, we can see this illustrated.

\begin{example}\label{ex:R+S+Ex3.5}
Consider the Riordan arrays $\Pas = \mathcal{R}\left(\frac{1}{1-t}, \frac{t}{1-t}\right)$ and $S = \mathcal{R}\left(\frac{1}{1-2t}, \frac{t}{1-t}\right)$ from Example \ref{ex:addRiordanAZsequence}. 
Recall that 

$$\Pas+S=\mathcal{R}\left(\frac{2-3t}{1-3t+2t^2}, \frac{t}{1-t}\right)=\begin{bmatrix}
2 & 0 & 0 & 0 & 0 & \cdots \\
3 & 2 & 0 & 0 & 0 & \cdots \\
5 & 5 & 2 & 0 & 0 & \cdots \\
9 & 10 & 7 & 2 & 0 & \cdots \\
17 & 19 & 17 & 9 & 2 & \cdots \\
\vdots & \vdots & \vdots & \vdots & \vdots & \ddots\\
\end{bmatrix}.
$$

Under the Riordan sumray definition, $\Pas+S=\mathcal{R}\left(\frac{1}{1-t}, \frac{1}{1-2t}, \frac{t}{1-t}, \frac{t}{1-t}\right)$.  These have shared $A$-sequence $(1,1,0,0,\dots)$.  Then
$$
\tilde{A}^2=\begin{bmatrix}
1 & 0 & 0 & 0 & \cdots\\
2 & 1 & 0 & 0 & \cdots\\
1 & 2 & 1 & 0 & \cdots\\
0 & 1 & 2 & 1 & \cdots\\
\vdots & \vdots & \vdots & \vdots & \ddots
\end{bmatrix}.
$$
Thus
$$
(\Pas+S)\tilde{A}^2=\begin{bmatrix}
2 & 0 & 0 & 0 & \cdots\\
7 & 2 & 0 & 0 & \cdots\\
17 & 9 & 2 & 0 & \cdots\\
36 & 26 & 11 & 2 & \cdots\\
\vdots & \vdots & \vdots & \vdots & \ddots
\end{bmatrix}=(\Pas+S)^{\sigma \sigma}.
$$
\end{example}

\section{The operations Der and Flip}
\label{sec:derflip} 

In this section we present two new operations on Riordan arrays. 
These are the {\em derivative} and the {\em flip}, which we define as follows. 

\begin{definition}
Let $\R(d(t),h(t))$ be a Riordan array. The {\em derivative} and {\em flip} are the following two operations respectively.
\begin{enumerate}
    \item $\Der$: $\R(d(t), h(t)) \mapsto \R(h'(t), t d(t));$
   \item $\Flip$: $\R(d(t), h(t)) \mapsto \R(h(t)/t, t d(t))$.  
\end{enumerate}
\end{definition}

The following observations are immediate. 

\begin{lemma} 
The operations $\Der$ and $\Flip$ satisfy the following. 
\begin{itemize}
    \item $\Flip$ is an involution.
    \item $\Flip(\Der(R))=\R(d(t), th'(t)).$
    \item $\Der(\Flip(R))=\R(d(t)+td'(t), h(t)).$
    \item Both $\Der$ and $\Flip$ map the Appell subgroup to the Lagrange subgroup, since \[\R(d(t), t)\overset{\Der,\Flip}{\longrightarrow}\R(1, td(t)).\]
    \item Both $\Der$ and $\Flip$ map the Lagrange subgroup to the Appell subgroup, 
    since 
    \[\Der(\R(1,h(t)))=\R(h'(t),t), \ \Flip(\R(1,h(t)))=\R(h(t)/t, t).\]
    \item $\Der$ maps the Derivative subgroup to the Bell subgroup, since 
     \[\R(h'(t),  h(t))\overset{\Der}{\longrightarrow} \R(h'(t), t h'(t)).\]
    \item $\Flip$ fixes an array $\R(d(t), h(t))$ if and only if it is in the Bell subgroup (when $h(t)=td(t)$), and here $\Der$ is the map  
    \[\R(d(t), t d(t))\overset{\Der}{\longrightarrow} \R(d(t)+td'(t), t d(t)).\]
    \item $\Der$ fixes an array $\R(d(t), h(t))$ if and only if $h'(t)=d(t)=h(t)/t$, i.e., if and only if $h(t)=ct$ for some nonzero constant $c$.
    \item $\Der(R)=\Flip(R)$ if and only if $h(t)=ct$ for some nonzero constant $c$, in which case both operations fix $R$, and $R=\R(c,ct)$.
\end{itemize}
\qed
\end{lemma}

From the definitions of $\Der$ and $\Flip,$ we have 
\begin{proposition}\label{prop:nk-entriesDerFlip} Let $R=\R(d(t), h(t))$ be any Riordan array.  The $(n,k)$-entries of $\Der$ and $\Flip$ are 
\[\Der(R)_{n,k}=[t^{n-k}] h'(t) d(t)^k, \quad \Flip(R)_{n,k}=[t^{n-k+1}] 
h(t) d(t)^k.\]
\qed
\end{proposition}

For an array $d_{n,k}$, the \emph{row sum} of the $n$th row is the sum $\sum_{k\ge 0} d_{n,k}$, while the \emph{alternating row sum} of the $n$th row is the signed sum 
$\sum_{\ge 0} (-1)^k d_{n,k}$. 

\begin{lemma}\label{lem:rowsums} 
Let $R=\R(d(t), h(t))$ be a Riordan array.
  The generating function for the row sums (respectively, alternating row sums) of 
\begin{enumerate}
    \item $\Der(R)$ is $\dfrac{h'(t)}{1-td(t)}$ (respectively, $\dfrac{h'(t)}{1+td(t)}$);
    \item $\Flip(R)$ is $\dfrac{h(t)}{t(1-td(t))}$ (respectively, $\dfrac{h(t)}{t(1+td(t))}$).
\end{enumerate}

\end{lemma}
\begin{proof}  This follows from the fact that the generating function for the row sums (respectively, alternating row sums)  of any Riordan array $R=\R(d(t), h(t))$ is $d(t)/(1-h(t))$ (respectively, $d(t)/(1+h(t))$), a consequence of Theorem~\ref{thm:FTRA} using the generating functions $f(t)=1/(1-t)$ (respectively, $f(t)=1/(1+t)$).
\end{proof}

The next two lemmas are simply applications of Theorem~\ref{thm:He-Sprugnoli}. As before, the bar denotes compositional inverse.
\begin{lemma}
\label{lem:derAseq}
If $R = \mathcal{R}(d(t), h(t))$ is a Riordan array, then the $A$-sequences of $\Der(R)$ and $\Flip(R)$ are given by
$$ A(t)=\dfrac{t}{\overline{t\textbf{d}(t)}},$$
and the $A$ sequence of $\Der^2(R)$ is 
$$ A(t) = \dfrac{t}{\overline{t\textbf{h}'(t) }}.$$
\end{lemma}


\begin{lemma}
\label{lem:derZseq}
If $R = \mathcal{R}(d(t), h(t))$ is a Riordan array, then the $Z$-sequences of $\Der(R)$ and $\Flip(R)$ are given by
$$ Z(t)=\dfrac{h'(\overline{t\textbf{d}(t)})-h'(0)}{ \overline{t\textbf{d}(t)}\cdot h'(\overline{t\textbf{d}(t)})}.$$
\end{lemma}

The following gives the connection between sums, derivatives, and flips of Riordan arrays.   Let $\Der^k$ (respectively, $\Flip^k$)  denote the operation $\Der$ (respectively, $\Flip$) iterated $k$ times. 

\begin{proposition}\label{prop:DerFlipSums}
Let $R_1=\R(d_1(t),h(t))$ and $R_2=\R(d_2(t),h(t))$ be Riordan arrays. Then the following hold:
\begin{itemize}
    \item[1.] $\Der^{2m}(R_1 + R_2) = \Der^{2m}(R_1) + \Der^{2m}(R_2),\  m\ge 1$, 
    \item[2.] $\Flip^{2m}(R_1 + R_2) = \Flip^{2m}(R_1) + \Flip^{2m}(R_2) ,\ m\ge 1$,
    \item[3.] $\Der(\Flip(R_1+R_2)) = \Der(\Flip(R_1)) + \Der(\Flip(R_2))$, and 
    \item[4.] $\Flip(\Der(R_1+R_2)) = \Flip(\Der(R_1)) + \Flip(\Der(R_1))$.
\end{itemize}    
\end{proposition}

\begin{proof}
These are easily checked. Note that for the first two items it suffices to check the case $n=1$.
\end{proof}

As noted before, $\Der$ maps the Appell subgroup to the Lagrange subgroup and vice versa. Therefore $\Der^{2m}$ maps the Appell subgroup to itself for any $m$.  
We have the following. 

\begin{theorem}\label{thm:derAppell} 
For a Riordan array $\R(d(t),t)$ in the Appell subgroup where $d(t) = d_0 + d_1t + d_2t^2 + \cdots $ and $m\geq 0$ we have 
$$ \Der^{2m}(R) = \R\left(\sum_{i= 0}^m S_{m+1,i+1} t^i d^{(i)}(t),t\right) = \R\left(\sum_{i\geq 0} d_i(i+1)^mt^i, t \right),$$
$$ \Der^{2m+1}(R) = \R\left(1, \sum_{i= 0}^m S_{m+1,i+1} t^{i+1} d^{(i)}(t)\right) = \R\left(1, t\sum_{i\geq 0} d_i(i+1)^mt^i \right),$$
where $S_{m+1,i+1}$ is the Stirling number of the second kind, and $d^{(i)}(t)$ denotes the $i$th derivative of $d(t)$. 
\end{theorem}
\begin{proof}
For each $m\geq 0$, the $\Der^{2m+1}(R)$ case follows immediately from $\Der^{2m}(R)$, so it suffices to show the latter. 
We begin by showing the first equality. 
Observe that when $m=0$, we have $\Der^{2m}(\R(d(t),t)) = \R(d(t), t)$. We then proceed by induction, assuming that $$\Der^{2(m-1)}(R) = \R\left(\sum_{i\geq 0}^{m-1} S_{m,i+1} t^i d^{(i)}(t),t\right).$$

Note that by the definition of $\Der$, we have
$$\Der^{2m}(R) = \R\left(\left(t\cdot \sum_{i\geq 0}^{m-1} S_{m,i+1} t^i d^{(i)}(t)\right)^{(1)},t\right).$$

It remains to simplify as follows.
\begin{align*}
    \left(t\cdot \sum_{i= 0}^{m-1} S_{m,i+1} t^i d^{(i)}(t)\right)^{(1)} &= \sum_{i= 0}^{m-1} [S_{m,i+1}(i+1) t^id^{(i)}(t) + S_{m,i+1}t^{i+1}d^{(i+1)}(t) ] \\
    &= \sum_{i= 0}^{m-1}S_{m,i+1}(i+1) t^id^{(i)}(t) + \sum_{i=1}^m S_{m,i}t^{i}d^{(i)}(t)\\
    &= S_{m,1}d(t) + S_{m,m}d^{(m)}(t) + \sum_{i=1}^{m-1} [S_{m,i+1}(i+1) + S_{m,i}]t^id^{(i)}(t) \\
    &=\sum_{i=0}^m S_{m+1,i+1} t^i d^{(i)}(t)
\end{align*}
The last equality uses the Stirling recurrence $S_{m+1,i+1} = (i+1)S_{m,i+1} + S_{m,i}$, and the fact that $S_{m,1} = S_{m+1,1} = S_{m,m}=S_{m+1,m+1}  = 1$. 

The second equality follows from the observation that 
\begin{align*}
    \sum_{j=0}^m S_{m+1,j+1} t^j d^{(j)}(t) &= \sum_{j=0}^m S_{m+1,j+1} \sum_{i\geq 0} d_i \dfrac{i!}{(i-j)!}t^i\\
    &= \sum_{i\geq 0} d_i \left(\sum_{j=0}^m S_{m+1,j+1}\dfrac{i!}{(i-j)!} \right) t^i \\
    &= \sum_{i\geq 0} d_i(i+1)^m t^i.
\end{align*}
For the final equality above, take the well-known identity  which counts functions from a set of size $(m+1)$  to a set of size $(i+1)$ according to the size $(j+1)$ of the image (see e.g. \cite[Eqn. (1.94d)]{StanEC1}), namely \[(i+1)^{m+1} = \sum_{j=0}^m (j+1)! \binom{i+1}{j+1}  S_{m+1,j+1} =\sum_{j=0}^m \frac{(i+1)!}{(i-j)!}  S_{m+1,j+1}, \]
and divide throughout by $(i+1)$. 
\end{proof}

\begin{corollary}\label{cor:nk-entriesAppell}
For any Riordan array $R=\R(d(t),t)$ in the Appell subgroup and $m\geq 0$, we have 
$$ \Der^{2m}(R)_{n,k} = d_{n-k}(n-k+1)^m.$$
\end{corollary}

\begin{corollary}\label{cor:b_n}
If $\Der^{2m}(\R(d(t),t))$ maps the sequence $(a_n)_{n\geq 0}$ to the sequence $(b_n)_{n\geq 0}$, then 
$$b_{n} =  \sum_{k\geq 0} a_kd_{n-k} (n-k+1)^m.$$
\end{corollary}

The operation $\Der$ can be given a combinatorial interpretation when applied to an element in the Appell subgroup. 
\begin{theorem}
\label{thm:DerAsComps}
Let $d(t) = d_0 + d_1t + d_2t^2 + \cdots$. Then $\Der (\R(d(t),t))_{n,k}$ is the number of \emph{weighted} compositions of $n$ with $k$ parts, where part $i$ has weight $d_{i-1}$. 
\end{theorem}

\begin{proof}
The following computation yields the result. 
\begin{align*}\Der(\R(d(t),t))_{n,k} &= \R(1,td(t))_{n,k} \\
&= [t^{n}] t^kd(t)^k \\
&= [t^{n}] (d_0t + d_1t^2+d_2t^3 + \cdots)^k \\
&= [t^{n}] \sum_{m\geq 0} \sum_{\substack{c_1+c_2+\cdots +c_k = m \\ c_1, c_2, \ldots, c_k \in \mathbb{Z}_{+}}} d_{c_1-1}d_{c_2-1}\cdots d_{c_k-1}t^{c_1}t^{c_2}\cdots t^{c_k} \\
&= \sum_{\substack{c_1+c_2+\cdots +c_k = n \\ c_1, c_2, \ldots, c_k \in \mathbb{Z}_{+}}} d_{c_1-1}d_{c_2-1}\cdots d_{c_k-1}.
\end{align*}
\end{proof}

\section{Applications}
\label{sec:applications}

In this section we investigate applications of the results in Section~\ref{sec:derflip} to well-known Riordan arrays. 
In particular, we apply $\Der$ and $\Flip$ to the Fibonacci array, Pascal array, Catalan array, and Shapiro array. In the process, we obtain various combinatorial identities.

\subsection{The Fibonacci Array}
We begin with a case in the Appell subgroup. 
Let $T=\R(1+t,t)$, and define $\Fib\coloneqq \Der(T) = \R(1, t+t^2)$. 
The Riordan array $\Fib$ is known as the Fibonacci array, as its row sums give the Fibonacci sequence. 
Indeed, using Lemma~\ref{lem:rowsums}, we see that the row sums of $\Fib$ have generating function $1/(1-t-t^2)$, which is the generating function for the Fibonacci sequence. 
The $(n,k)$-entry in $\Fib$ is given by 
\begin{equation}
    \Fib_{n,k}=[t^{n-k}] (1+t)^k= \binom{k}{n-k}, n\ge k\ge 0.
\end{equation}
$$
\,\,\,\,\qquad T=\begin{bmatrix}
1 & 0 & 0 & 0 & 0 & 0 & \cdots\\
1 & 1 & 0 & 0 & 0 & 0 & \cdots\\
0 & 1 & 1 & 0 & 0 & 0 & \cdots\\
0 & 0 & 1 & 1 & 0 & 0 & \cdots\\
0 & 0 & 0 & 1 & 1 & 0 & \cdots\\
0 & 0 & 0 & 0 & 1 & 1 & \cdots\\
\vdots & \vdots & \vdots & \vdots & \vdots & \vdots & \ddots
\end{bmatrix},
\mathrm{Fib}=\R(1, t+t^2)=\begin{bmatrix}
1 & 0 & 0 & 0 & 0 & 0 & \cdots\\
0 & 1 & 0 & 0 & 0 & 0 & \cdots\\
0 & 1 & 1 & 0 & 0 & 0 & \cdots\\
0 & 0 & 2 & 1 & 0 & 0 & \cdots\\
0 & 0 & 1 & 3 & 1 & 0 & \cdots\\
0 & 0 & 0 & 3 & 4 & 1 & \cdots\\
\vdots & \vdots & \vdots & \vdots & \vdots & \vdots & \ddots
\end{bmatrix}
$$
Furthermore, $\Der^2(T) =\R(1+2t, t)$ and $\Der^3(T)=\R(1, t+2t^2)$ are the following Riordan arrays.
$$
\Der^2(T)\!=\!\begin{bmatrix}
1 & 0 & 0 & 0 & 0 & 0 & \cdots\\
2 & 1 & 0 & 0 & 0 & 0 & \cdots\\
0 & 2 & 1 & 0 & 0 & 0 & \cdots\\
0 & 0 & 2 & 1 & 0 & 0 & \cdots\\
0 & 0 & 0 & 2 & 1 & 0 & \cdots\\
0 & 0 & 0 & 0 & 2 & 1 & \cdots\\
\vdots & \vdots & \vdots & \vdots & \vdots & \vdots & \ddots
\end{bmatrix}\!,
\, \mathrm{Jac} :=\, \Der^3(T)\! =\!\begin{bmatrix}
1 & 0 & 0 & 0 & 0 & 0 & \cdots\\
0 & 1 & 0 & 0 & 0 & 0 & \cdots\\
0 & 2 & 1 & 0 & 0 & 0 & \cdots\\
0 & 0 & 4 & 1 & 0 & 0 & \cdots\\
0 & 0 & 4 & 6 & 1 & 0 & \cdots\\
0 & 0 & 0 & 12 & 8 & 1 & \cdots\\
\vdots & \vdots & \vdots & \vdots & \vdots & \vdots & \ddots
\end{bmatrix}.
$$
Note that by Theorem~\ref{thm:FTRA}, $\Der^2(T)$ has row sums $(1,3,3,\ldots)$ with generating function $(1+2t)/(1-t)$, while 
the row sums of $\Der^3(T)$ have generating function $1/(1-t-2t^2)$. The latter is the generating function for the Jacobsthal numbers, which are recursively defined by $a_n = a_{n-1} + 2a_{n-2}$ with initial values $a_0=a_1=1$ (see \cite[A001045]{OEIS}). 
For this reason we define $\Jac :=\Der^3(T)$ to be the Jacobsthal array.

We also observe that 
\begin{itemize}
    \item $\Flip(\Der(T))) = \Flip(\Fib) = \Flip(\R(1,t+t^2)) = T$, and 
    \item $\Flip(\Der^3(T)) = \Flip(\mathrm{Jac}) = \Flip(\R(1,t+2t^2)) = \Der^2(T)$.
\end{itemize}

Since $\mathrm{Jac}=\R(1, t(1+2t)),$ we have that $\mathrm{Jac}_{n,k}=[t^{n-k}](1+2t)^k=2^{n-k}\binom{k}{n-k}.$
The $(n+1)$-st Jacobsthal number is thus $\sum_{k=0}^n 2^{n-k}\binom{k}{n-k}$. 
On the other hand,  the generating function $1/(1-t-2t^2)=1/((1+t)(1-2t))$ gives the explicit formula $J_n=\frac{1}{3}(2^{n+1}+(-1)^n)$ 
for the $n$th Jacobsthal number. 
These identities are given in \cite[A001045]{OEIS}.  

We give a bijective proof of the following combinatorial interpretation of the Jacobsthal numbers  in terms of pattern avoidance. 

\begin{theorem}
Let $D_{n}(231,132)$ be the set of derangements of length $n$ avoiding the patterns $231$ and $132$. If $J_{n}$ denotes the $n$th Jacobsthal number, then 
$$
J_{n+1} = \sum_{k=0}^{n}2^{n-k}\binom{k}{n-k} = |D_{n+2}(231,132)|.
$$
\end{theorem}
\begin{proof}(Combinatorial)
First, note that every $\pi\in D_{n+2}(231,132)$ must begin with $n+2$ since if we have elements both before and after $n+2$, then a forbidden pattern will be formed. Furthermore, $n+2$ cannot be at the last position, as $\pi$ is a derangement. In addition, to avoid $231$ and $132$, all the numbers preceding $1$ in $\pi$ must be in decreasing order and all the numbers following $1$ in $\pi$ must be in increasing order. For example, $53124\in D_{5}(231,132)$. This implies that there is a unique index $i\geq 2$, such that $\pi_{i}<i$ and $\pi_{j}>j$ for $j = 1, \ldots , i-1$. For a fixed $i$, we have that the numbers $\pi_{2},\ldots , \pi_{i-1}$ are $i-2$ numbers in decreasing order among $i, i+1 , \ldots , n+1$. Thus, we can select these $i-2$ numbers and determine the segment $\pi_{2}\cdots\pi_{i-1}$ in $\binom{(n+1)-(i-1)}{i-2} = \binom{n+2-i}{i-2}$ ways. The non-selected numbers in $\{i, i+1 , \ldots , n+1\}$ must be in increasing order at the end of $\pi$. What remains is to determine the positions of $i-1,\ldots , 2, 1$. Note that if we do this sequentially for each of the listed numbers, every time we will have exactly two choices for that position - the leftmost or the rightmost unoccupied position in the permutation. The only exception is the position of $1$ for which we will have only one possible choice at the end. Hence we shall multiply by $2^{i-2}$. For instance, when $n=3$ and $i=3$, we will have $\binom{5-3}{1}2^{1}$ such permutations in $D_{5}(231,132)$: $54123, 54213, 53124, 53214$. The number at position $2$, i.e., $\pi_{2}$ is determined in $\binom{5-3}{1} = 2$ ways since it can be $3$ or $4$. If, for example, $\pi_{2} = 3$, then the number $4$ must be at the last position, i.e., $\pi_{5}=4$. Then, we have $2$ choices for the position of the number $2$ - either after $3$ or before $4$. The number $1$ must be at the last unoccupied position. Note that by following this simple algorithm for construction of $\pi$, we always obtain a derangement. Summing over the possible values of $i$, we get 
$$
|D_{n+2}(231,132)| = \sum\limits_{i=2}^{n+2} \binom{n+2-i}{i-2}2^{i-2} = \sum_{j=0}^{n} \binom{n-j}{j}2^{j} = \sum_{k=0}^{n}2^{n-k}\binom{k}{n-k},
$$
as claimed.
\end{proof}

Applying Theorem~\ref{thm:DerAsComps} to $\R(1+ct,t)$ shows that the previous identity is a special case of the following. 

\begin{identity}
\label{id:12comps}
Let $a_n$ denote the number of compositions of $n$ using parts $1$ and $2$, with $c$ available colors for the $2$'s. Then
\[a_n = \sum_{k\geq 0}  c^{n-k}\binom{k}{n-k}. \]
\end{identity}

\begin{proof}
There are $\binom{k}{n-k}$ compositions of $n$ with $k$ parts from $\{1,2\}$ of which $n-k$ are $2$'s. We have $c^{n-k}$ ways to color these $2$'s. 
Summing over the possible number of parts gives the result.
\end{proof}

\subsection{The Pascal Array} 

The Pascal array $\Pas$ and its inverse $\Pas^{-1}$ are the Riordan arrays 
\begin{equation*}\label{eqn:PascalArray} \Pas=\R\left(\frac{1}{1-t}, \frac{t}{1-t}\right) \qquad \text{ and } \qquad 
\Pas^{-1}=\R\left(\frac{1}{1+t}, \frac{t}{1+t}\right).
\end{equation*}
Their general terms are given by $\Pas_{n,k}=\binom{n}{k}$ and $\Pas^{-1}_{n,k}=(-1)^{n-k}\binom{n}{k}$ respectively.

Recall from \cite[Sec.~1.3]{Petersen2015} that  the \emph{Eulerian number} $A(n,k)$  counts the number of permutations on $n$ letters with exactly $k$  descents, where $i$ is a descent of a permutation $\sigma$ on $n$ letters if and only if $\sigma(i)>\sigma(i+1), 1\le i\le n-1.$ Denote by $des(\sigma)$ the number of descents of the permutation $\sigma$. The generating function for the Eulerian numbers $\{A(n,k)\}_{k=0}^{n-1}$ is the \emph{Eulerian polynomial} $A_n(t):=\sum_{\sigma\in S_n} t^{des(\sigma)}=\sum_{k=0}^{n-1} A(n,k)t^k, n\ge 1.$ We define $A_0(t):=1.$  Thus $A_1(t)=1, A_2(t)=1+t,$ 
$A_3(t)=1+4t+t^2, A_4(t)=1+11t+11t^2+t^3.$
\footnote{Note: Stanley \cite[Sec.~1.4]{StanEC1} and others define the Eulerian polynomial  to be $\sum_{\sigma\in S_n} t^{1+des(\sigma)}=t\cdot A_n(t)$.}

We have the following recurrence \cite[Theorem 1.4]{Petersen2015} for the Eulerian polynomials:
\begin{equation}\label{eqn:EulerPolyRec}
A_{n+1}(t)=(1+nt)A_n(t)+t(1-t)A_n'(t).
\end{equation}

\begin{theorem}\label{thm:nk-entriesDerPas}
Let $A_i(t)$ denote the $i$th Eulerian polynomial. For each $i\geq 0$, we have 
\begin{align*}
\Der^{2i}(\Pas) &= \R\left( \dfrac{A_i(t)}{(1-t)^{i+1}}, \dfrac{tA_i(t)}{(1-t)^{i+1}} \right),  \\
\Der^{2i+1}(\Pas) &= \R\left( \dfrac{A_{i+1}(t)}{(1-t)^{i+2}}, \dfrac{tA_i(t)}{(1-t)^{i+1}} \right), \\
\Der^{2i}(\Pas^{-1}) &= \R\left( \dfrac{A_i(-t)}{(1+t)^{i+1}}, \dfrac{tA_i(-t)}{(1+t)^{i+1}} \right),  \\
\Der^{2i+1}(\Pas^{-1}) &= \R\left( \dfrac{A_{i+1}(-t)}{(1+t)^{i+2}}, \dfrac{tA_i(-t)}{(1+t)^{i+1}} \right). 
\end{align*}
\end{theorem}

\begin{proof}
Using the well-known generating function \cite[Corollary~1.1]{Petersen2015}
\begin{equation}\label{eqn:powers-of-integers-EulerianPoly}\sum_{k=0}^\infty k^i t^k=\frac{t A_i(t)}{(1-t)^{i+1}}
\end{equation}
we observe that 
$$ \dfrac{d}{dt}\left(\dfrac{tA_i(t)}{(1-t)^{i+1}} \right) = \sum_{k=1}^\infty k^{i+1} t^{k-1}= \dfrac{A_{i+1}(t)}{(1-t)^{i+2}}.$$
The equations for the derivatives of $\Pas$ and $\Pas^{-1}$ then follow readily by induction on $i$.  
\end{proof}

\begin{example}\label{ex:DerFlipPascal}
$$\Der(\Pas)=\R\left(\dfrac{1}{(1-t)^2}, \dfrac{t}{1-t}\right) = \begin{bmatrix}
1 & 0 & 0 & 0 & 0 &  \cdots\\
2 & 1 & 0 & 0 & 0 &  \cdots\\
3 & 3 & 1 & 0 & 0 &  \cdots\\
4 & 6 & 4 & 1 & 0 &  \cdots\\
5 & 10 & 10 & 5 & 1 &  \cdots\\
\vdots & \vdots & \vdots & \vdots & \vdots & \ddots
\end{bmatrix},$$
$$\Der^2(\Pas)=\R\left(\dfrac{1}{(1-t)^2}, \dfrac{t}{(1-t)^2}\right) = 
\begin{bmatrix}
1 & 0 & 0 & 0 & 0 &  \cdots\\
2 & 1 & 0 & 0 & 0 &  \cdots\\
3 & 4 & 1 & 0 & 0 &  \cdots\\
4 & 10 & 6 & 1 & 0 &  \cdots\\
5 & 20 & 21 & 8 & 1 & \cdots \\
\vdots & \vdots & \vdots & \vdots & \vdots & \ddots
\end{bmatrix}, 
$$
$$\Flip(\Der(\Pas))=\R\left(\dfrac{1}{1-t}, \dfrac{t}{(1-t)^2}\right) = 
\begin{bmatrix}
1 & 0 & 0 & 0 & 0 &  \cdots\\
1 & 1 & 0 & 0 & 0 &  \cdots\\
1 & 3 & 1 & 0 & 0 &  \cdots\\
1 & 6 & 5 & 1 & 0 &  \cdots\\
1 & 10 & 15 & 7 & 1 & \cdots \\
\vdots & \vdots & \vdots & \vdots & \vdots & \ddots
\end{bmatrix}.
$$

We can compute the general terms of these Riordan arrays to be
\[\Der(\Pas)_{n,k}=\binom{n+1}{k+1}, \,\, \Der^2(\Pas)_{n,k}=\binom{n+k+1}{n-k}, \,\, \Flip(\Der(\Pas))_{n,k} = \binom{n+k}{2k}.\]
Using Lemma~\ref{lem:rowsums}, the row sums of $\Der^2(\Pas)$ and $\Flip(\Der(\Pas))$ have the  generating functions
$$ \dfrac{1}{1-3t+t^2} \quad \text{ and } \quad \dfrac{1-t}{1-3t+t^2}, $$
respectively. These are well known to be the generating functions for the two bisections of the Fibonacci sequence, $\sum_{n\ge 0} F_{2n+1} t^n$, i.e. $(1,3,8,21,...)$, and $\sum_{n\ge 0} F_{2n} t^n$, i.e. $(1,2,5,13,...)$.
See also Example~\ref{ex:Fib-bisection}.

Similarly, the alternating row sums have generating functions as follows:

For $\Pas$ itself, the alternating row sum generating function is 1, reflecting the fact that the alternating sum of the binomial coefficients in the $n$th row of Pascal's triangle is 0 for $n\ge 1$. 
Note that the 
alternating row sums of $\Pas$ are the row sums of $\Pas^{-1}$.

For $\Der(\Pas)$, it  is $1/(1-t)$, i.e., the alternating sum in each row is 1.

For $\Der^2(\Pas)$, it is $(1+t)/(1+t^3)$, and hence the alternating sum along the $n$th row is $\sum_{n\ge 0} (t^{3n}+t^{3n+1}),$ i.e. it equals 
\[\begin{cases} (-1)^{\lfloor\frac{m}{3}\rfloor}, & \text{if $m\equiv 0,1 \hspace{-0.3cm} \mod 3$}, \\
0, &\text{otherwise.}
\end{cases}
\]

For $\Flip\Der(\Pas))$, it is $(1-t^2)/(1+t^3)$, and hence the alternating sum along the $n$th row equals 
\[\begin{cases} (-1)^{\frac{m}{3}}, & \text{if $m$ is divisible by 3,}\\
(-1)^{\frac{m+1}{3}}, & \text{if $m+1$ is divisible by $3$,}\\
0, & \text{otherwise.}
\end{cases}\]

\end{example}

\begin{definition}\label{def:INVERT}\cite{BernsteinSloane1995}
Define the INVERT transform  of the sequence $(a_n)_{n\ge 1}$ to be the sequence $(b_n)_{n\ge 1}$ where 
\[1+\sum_{n\ge 1} b_n t^n=(1-\sum_{n\ge 1} a_n t^n)^{-1}.\]
\end{definition}

We have the following interesting relationship between the even derivatives of the Riordan array $\Pas$ and the INVERT transform of the sequence of $n$th powers of the positive integers, for  fixed $n$.

\begin{proposition}\label{prop:Der-Pas-INVERT-transformS}
Let $n$ be a fixed positive integer.  Suppose the sequence $(b_k)_{k\ge 1}$ is the INVERT transform of the sequence $(k^n)_{k\ge 1}$.  Then the generating function for the row sums of $\Der^{2n}(\Pas)$ is 
\[\sum_{k\ge 0} b_{k+1} t^k.\]
\end{proposition}
\begin{proof} By definition, we have 
\[\frac{1}{1-\sum_{k\ge 1} k^n t^k}=1+\sum_{k\ge 1} b_k t^k.\]
From Theorem~\ref{thm:nk-entriesDerPas}, we have 
\[\Der^{2n}(\Pas)=\R\left(\frac{A_n(t)}{(1-t)^{n+1}}, \frac{tA_n(t)}{(1-t)^{n+1}}\right)
=\R\left(\sum_{k\ge 1} k^n t^{k-1}, \sum_{k\ge 1} k^n t^{k}\right).\]
Write $f_n(t)=\sum_{k\ge 1} k^n t^{k}.$  
The generating function for row sums is thus 
\begin{equation*}
    \frac{t^{-1}f_n(t)}{1-f_n(t)}
    =t^{-1}\left(\frac{1}{1-f_n(t)}-1\right)
    =t^{-1}\left(\sum_{k\ge 1} b_k t^k\right),
\end{equation*}
which equals $\sum_{j\ge 0} b_{j+1} t^j$, as claimed.
\end{proof}

\begin{example}\label{ex:Der-Pas-row-sumS}
The generating functions of the first few row sums of $\Der^{m}(\Pas)$ for even values of  $m$, together with the relevant sequences in OEIS, are as follows:

$\Der^2(\Pas):$  row sum generating function $\dfrac{1}{1-3t+2t^2}$. (\cite[A000225]{OEIS})

$\Der^4(\Pas):$  row sum generating function $\dfrac{1+t}{1-4t+2t^2-t^3}$.
(\cite[A033453]{OEIS})

$\Der^6(\Pas):$  row sum generating function $\dfrac{1+4t+t^2}{1-5t+2t^2-5t^3+t^4}$. (\cite[A144109]{OEIS})

Many other connections appear to hold.  For instance, 
the generating function for the row sums of $\Der^3(\Pas)$ is 
\[\dfrac{1+t}{(1-t)(1-3t+t^2)}=\dfrac{1+t}{1-4t+4t^2-t^3},\]
which is two less than the bisection of Lucas numbers \cite[A004146]{OEIS}.

The row sums of $\Der^2(\Pas) + \Der^3(\Pas)$ have generating function $2/(1-4t+4t^2-t^3)$, and are twice the sequence \cite[A027941]{OEIS} of the odd Fibonacci minus one:  $(F_{2n + 1} - 1)_{n=0}^\infty$. 

More generally, one can add the $(2i)$th and $(2i+1)$th derivatives, since they have the same $h(t)$.  The sum is the Riordan array
$\mathcal{R}( \sum_{\ell\ge 0} ((\ell+1)^i+\ell^{i+1}) t^\ell, \sum_{\ell\ge 1} \ell^i t^\ell),$ and the row sum generating function in terms of the Eulerian polynomials is \[\frac{(1-t)A_i(t)+A_{i+1}(t)}{(1-t)^{i+2}-tA_i(t)}.\]
A similar formula can be computed for the row sum generating function of the sum 

\noindent
$\Flip(\Der^{2n-1}(\Pas))+\Der^{2n}(\Pas).$
\end{example}

\subsection{The Catalan Array}\label{subsec:Catalan}

\cite{StanCatalan2015}
Let $C(t)$ denote the generating function $\sum_{n\ge 0} C_n t^n$ for the Catalan sequence $(1,1,2,5,14,\ldots)$. The defining equation for $C(t)$ is 
\[t C^2(t) -C(t) +1=0,\]
giving the well-known formula 
\[C(t)=\frac{1-\sqrt{1-4t}}{2t}, \text{ and hence also } C(t)-1=\frac{1-2t-\sqrt{1-4t}}{2t}.\]
 Recall that  the Catalan array $\Cat$ is the Riordan array defined by
\begin{equation}\label{eqn:CatalanArray} \Cat\coloneqq\R(C(t), t C(t)).
\end{equation}
Its inverse $\Cat^{-1}$ is \cite[Theorem~5.2]{LuzonMerliniMoronSprugnoli2012LAA}
\begin{equation}\label{eqn:CatalanArrayInv}
\Cat^{-1}=\R(1-t, t-t^2).
\end{equation}

From \cite[Theorems~5.2-5.3]{LuzonMerliniMoronSprugnoli2012LAA}, we have  $\Cat^{-1}_{n,k}=(-1)^{n-k}\binom{k+1}{n-k}$ and hence 
    \[ \Cat_{n,k}=\frac{k+1}{n+1}\binom{2n-k}{n-k}. \]

\begin{theorem}\label{thm:nk-entriesDerCat} 
For the Catalan array $\Cat=\R(C(t), tC(t))$ we have
\begin{equation*} \Der(\Cat)=\R\left(\frac{1}{\sqrt{1-4t}},\frac{1-\sqrt{1-4t}}{2}\right) \text{ and }
(\Der(\Cat))^{-1}=\R\left(1-2t, t-t^2\right), \end{equation*}
and hence 
\[\Der(\Cat)_{n,k}=\binom{2n-k}{n-k}.\]
\end{theorem}
\begin{proof}
We mimic the ingenious method of \cite{MERLINI2017160}. Using Theorem~\ref{thm:Shapiro}, we first calculate the inverse array $R=(\Der(\Cat))^{-1}$ of $\Der(\Cat)$. Then we apply Theorem~\ref{thm:InvRiordanElements} to compute  the $(n,k)$-entry of the inverse of 
$R$, which is of course precisely the $(n,k)$-entry of $\Der(\Cat)$.


This scheme exploits the fact that the inverse of $\Der(\Cat)$ involves only simple polynomials, and hence its  entries are easier to compute.  

From the defining recurrence $tC^2(t)=C(t)-1,$ we observe that 
\begin{enumerate}
\item[(0)] $t\,C(t) =\frac{1-\sqrt{1-4t}}{2}=\frac{2t}{1+\sqrt{1-4t}}$;
\item $\frac{d}{dt} (t\,C(t))= 
    \frac{1}{\sqrt{1-4t}}$;
    \item the compositional inverse of $t\,C(t)$ is $t-t^2$, and 
    \item $ \frac{d}{dt} (t\,C^2(t))=C'(t).$
\end{enumerate}    
Hence we have $\Der(\Cat)=\R(d^*(t), h^*(t))$, where $d^*(t)=(\sqrt{1-4t})^{-1}$ and $h^*(t)=t C(t)$, and therefore $\bmbar{h^*}(t)=t-t^2.$ 
Clearly \[\frac{1}{d^*(t-t^2)}=1-2t=\bmbar{h^*}'(t), \] and hence the inverse of $\Der(\Cat)$ is $\R(d(t), h(t))$ where $d(t)=1-2t$ and  $h(t)=t-t^2=\bmbar{h^*}(t).$
Now apply Theorem~\ref{thm:InvRiordanElements}, noting that $h(t)/t=1-t$ and $d(t)=h'(t)$ in this case.
Extracting the $(n,k)$-entry of $\Der(\Cat)$ can now be done using Definition~\ref{def:riordan}.
\end{proof}

\begin{example}\label{ex:AseqDerCat-FlipCat}
The compositional inverse of $t/(1-t)^2$ is $1-C(-t)$, and hence by Theorem~\ref{thm:He-Sprugnoli}, the generating function for the $A$-sequence of both $\Der(\Cat)$ and $\Flip(\Cat)$ is  
\begin{equation*}
\begin{split}\frac{1}{2}(1+2t+\sqrt{1+4t})&=\dfrac{t}{1-C(-t)}=1+t(1+C(-t))\\ 
&= 1+2t-t^2+2t^3-5 t^4+14t^5-42t^6+\cdots
\end{split}
\end{equation*}
\end{example}

Recall \cite{StanEC1} that a \emph{weak composition} of a nonnegative integer $m$ is a finite sequence of nonnegative integers whose sum is $m$. 
\begin{corollary} For each fixed $k\ge 0$, 
\begin{equation}\label{eqn:DerCat-cols1}
C^k(t)=\sqrt{1-4t}\sum_{n\ge k} \binom{2n-k}{n} t^{n-k},
\end{equation}
or equivalently
\begin{equation}\label{eqn:DerCat-cols2}
C^k(t)\sum_{n\ge 0}\binom{2n}{n}t^n=\sum_{n\ge k} \binom{2n-k}{n} t^{n-k}.
\end{equation}
See \cite[Exercise A.32 (a), (b)]{StanCatalan2015} for \eqref{eqn:DerCat-cols1} and also for the power series for $C^k(t)$.
This gives the following identity. 
\begin{identity}
\label{id:weakCompCatalan}
\[\sum_I \binom{2i_0}{i_0}C_{i_1} C_{i_2}\ldots C_{i_k}=\binom{2n-k}{n},\]
where the sum on the left runs over all weak compositions $I=(i_0,i_1,\ldots i_k)$ of $n-k$, for fixed $k\ge 0.$
\end{identity}
\end{corollary}




From Theorem~\ref{thm:nk-entriesDerCat} we know that $\Der(\Cat)_{n,k} = \binom{2n-k}{n-k}$. 
Lemma~\ref{lem:rowsums} gives, for the row sums, the generating function 
\[\frac{1-\sqrt{1-4t}}{2t\sqrt{1-4t}}=\frac{1}{2}\sum_{n\ge 0} \binom{2(n+1)}{n+1}t^n. \]
In particular the $n$th row sum of $\Der(\Cat)$ is $\frac{1}{2}\binom{2(n+1)}{n+1}=\frac{n+2}{2}\cdot C_{n+1}$.

Another Riordan array related to the Catalan numbers is the Shapiro array\cite{MERLINI2017160}  $\Sha=\R(C(t), t C^2(t))=\R(C(t), C(t)-1)$. 
We will consider the Shapiro array in more detail in Section~\ref{subsec:shapiro}, but for now we recall (\cite{MERLINI2017160}, \cite{Radoux2000JCompApplMath}) the following property of the Shapiro array:
\[\Sha\cdot(1,3,5,7,...)^t=(1,4,4^2,4^3,\ldots)^t.\]
See \cite[Theorem~2.1]{MERLINI2017160} for a combinatorial proof. 
We obtain an analogous result for $\Der(\Cat)$, as Theorem~\ref{thm:FTRA} now tells us that $\Der(\Cat)$ transforms powers of $2$ into powers of $4$:
\[\Der(\Cat)\cdot (1,2,2^2,2^3,\ldots)^t=(1,4,4^2,4^3,\ldots)^t.\]
Hence we have the following.

\begin{identity}
\label{id:powerOf4}
$$\sum_{k\geq 0} 2^k\binom{2n-k}{n-k} = 4^n= \sum_{k\ge 0} (2k+1) \binom{2n-k+2}{n-k}$$
\end{identity}
\begin{proof} (Combinatorial) 
We can rewrite the first equation as 
$$
2\sum\limits_{k=1}^{n} 2^{k-1}\binom{2n-k}{n} = 2^{2n} - \binom{2n}{n}.
$$
Then we can proceed as in the proof of Identity~\ref{id:derSH}.
\end{proof}

Similarly, since the sequence $(n+1)_{n=0}^\infty$ has generating function $(1-t)^{-2}$, Theorem~\ref{thm:FTRA} tells us that the generating function for $\Der(\Cat)\cdot (1,2,3,4,5,...)^t$ is 
\[ \frac{1}{\sqrt{1-4t}} \left(1- \frac{1-\sqrt{1-4t}}{2}\right)^{-2} \!=\!  \frac{(1-2t)+\sqrt{1-4t}}{2t^2\sqrt{1-4t}}\!=
\frac{1}{2t^2}\left((1-2t)\sum_{n\ge 0}\binom{2n}{n} t^n -1\right), \]
and hence the coefficient of $t^n, n\ge 0, $ is 
$\binom{2n+2}{n}.$  We therefore obtain:
\begin{identity}
\label{id:OddsToCat}
$$\sum_{k\geq 0} (k+1)\binom{2n-k}{n-k} = \binom{2(n+1)}{n} = (n+1)\cdot C_{n+1}$$
\end{identity}
\begin{proof} (Combinatorial) 
Rewrite the identity as 
$$
\sum\limits_{k=1}^{n} (k+1)\binom{2n-k}{n} = \binom{2n+2}{n+2} - \binom{2n}{n}.
$$
The right-hand side counts all of the subsets of $[2n+2]$ with exactly $n+2$ elements, which do not contain both $2n+1$ and $2n+2$. Thus, each such subset $A$ contains at least $n+1$ numbers in $[2n]$. Let the $(n+1)$th largest number in $A$ be $2n-k+1$, where $k$ can be between $1$ and $n$. The smallest $n$ numbers in $A$ can be chosen in $\binom{2n-k}{n}$ ways and the largest number in $A$ can be chosen in $(k+1)$ ways.
\end{proof}

Combining the expression for the row sums of $\Der(\Cat)$ with the preceding sum (simply using $2k+1=2(k+1)-1$) we can conclude that the $n$th term of 
 $\Der(\Cat)\cdot (1,3,5,7,9,...)^t$  equals  $\frac{3n+2}{n+2} \binom{2n+1}{n+1}$, which counts the number of positive clusters of Type $D_{n+2}$ \cite{FZ03}. 
Thus we obtain the following.
\begin{identity}
\label{id:pureDescents}
$$\sum_{k\geq 0} (2k+1)\binom{2n-k}{n-k} = \frac{3n+2}{n+2} \binom{2n+1}{n+1}$$
\end{identity}
\begin{proof} (Combinatorial)

Subtracting the previous Identity~\ref{id:OddsToCat} and using $\binom{2n+2}{n+2} = \frac{2n+2}{n+2}\binom{2n+1}{n+1}$, we see that it suffices to prove the following:
$$
\sum_{k\geq 1} k\binom{2n-k}{n}  = \sum_{k\geq 1} k\binom{2n-k}{n-k} =  \frac{n}{n+2} \binom{2n+1}{n+1} = \binom{2n+1}{n+2}.
$$
Now, let us consider a choice of $n+2$ numbers out of $2n+1$ numbers labeled with $1,2, \ldots , 2n+1$ and let $2n-k+1$ be the second largest label of a selected number. Note that $k\geq 1$. Out of the numbers with labels $1,2, \ldots 2n-k$, exactly $n$ must be selected and this can happen in $\binom{2n-k}{n}$ ways. The selected number with the largest label can be chosen in $k$ different ways since its label can be each of $2n-k+2, 2n-k+3, \ldots 2n+1$.
\end{proof}
The right-hand side is the sequence \cite[A129869]{OEIS}. Note that this is also the total number of all pure descents \cite{BarilKirgizov2017} 
whereas the right-hand side of Identity~\ref{id:derSH} equals the total number of inversions in all $321$-avoiding permutations of length $n$.

Next we apply $\Der$ again to $\Der(\Cat)$, obtaining
$$\Der^2(\Cat)=\R\left(\frac{1}{\sqrt{1-4t}}, \frac{t}{\sqrt{1-4t}}\right) = 
\begin{bmatrix}
1 & 0 & 0 & 0 & 0 & 0 & \cdots\\
2 & 1 & 0 & 0 & 0 & 0 & \cdots\\
6 & 4 & 1 & 0 & 0 & 0 & \cdots\\
20 & 16 & 6 & 1 & 0 & 0 & \cdots\\
70 & 64 & 30 & 8 & 1 & 0 & \cdots\\
252 & 256 & 140 & 48 & 10 & 1 & \cdots\\
\vdots & \vdots & \vdots & \vdots & \vdots & \vdots & \ddots
\end{bmatrix}.
$$ 

The Riordan array $\Der^2(\Cat)$ is also of combinatorial interest, as its $(n,k)$-entries count certain lattice paths as described by the following theorem. 
The number of such paths is given by \cite[A026671]{OEIS}.

\begin{theorem}
\label{thm:latticePaths}
Let $a_n$ denote the number of lattice paths from $(0,0)$ to $(n,n)$ with steps $(0,1)$, $(1,0)$ and, when on the diagonal, $(1,1)$.
If $a_{n,k}$ denotes the number of such paths with exactly $k$ diagonal steps, then 
$$ \Der^2(\Cat)_{n,k} = a_{n,k} = 4^{n-k}\binom{n-\frac{k+1}{2}}{n-k}$$
and
$$a_n = \sum_{k=0}^n 4^{n-k}\binom{n-\frac{k+1}{2}}{n-k}.$$
\end{theorem}

\begin{proof}
First note that by definition, we have 
\[\Der^2(\Cat)_{n,k} =[t^{n-k}]\left\lbrace \dfrac{1}{\sqrt{1-4t}} \cdot \left(\dfrac{1}{\sqrt{1-4t}}\right)^k \right\rbrace =[t^{n-k}](1-4t)^{-\frac{k+1}{2}},\]
and hence 
\[\Der^2(\Cat)_{n,k}=(-4)^{n-k}\binom{-\frac{k+1}{2}}{n-k}\]
which is as claimed, using the well-known formula  $\binom{-m}{j}=(-1)^j\binom{m+j-1}{j}$.

Now we observe that the first column of $\Der^2(\Cat)$ consists of the central binomial coefficients, as it has generating function $1/\sqrt{1-4t}$. The central binomial coefficients count the number of lattice paths from $(0,0)$ to $(n,n)$ using no diagonal steps, so $\Der^2(\Cat)_{n,0} = a_{n,0}$.
We fix $n$ and proceed by induction on $k$. Assume $\Der^2(\Cat)_{n,k-1} = a_{n,k-1}$. 
Then we have, for $k\geq 1$, 

\begin{align*}
\Der^2(\Cat)_{n,k} 
&= [t^{n}]\left\lbrace \dfrac{1}{\sqrt{1-4t}} \cdot \left(\dfrac{t}{\sqrt{1-4t}}\right)^k \right\rbrace  \\
&=[t^{n-1}] \left\lbrace \dfrac{1}{\sqrt{1-4t}} \cdot \left(\dfrac{1}{\sqrt{1-4t}}\cdot \left(\dfrac{t}{\sqrt{1-4t}} \right)^{k-1} \right)\right\rbrace  \\
&= \sum_{i=0}^{n-1} a_{i,0} \Der^2(\Cat)_{n-1-i,k-1} \\
&= \sum_{i=0}^{n-1} a_{i,0} a_{n-i-1,k-1}.    
\end{align*}
The set of lattice paths with $k$ steps of the form $(1,1)$ on the diagonal can be partitioned according to where their final diagonal step lies. 
Let $L_i$ denote the set of lattice paths whose first diagonal step is the step from $(i,i)$ to $(i+1,i+1)$. 
The lattice paths in $L_i$ then consist of the paths of the form $P_1 D P_2$, where $P_1$ is a path to $(i,i)$ with no diagonal steps, $D$ is the first diagonal step in the path, and $P_2$ is a path of length $n-i-1$ with $k-1$ steps on the diagonal. 
The paths in $L_i$ are thus counted by $a_{i,0} a_{n-i-1,k-1}$. 
Summing over all $i$ gives $a_{n,k}$, as desired.
\end{proof}

Applying $\Flip$ to $\Der(\Cat)$ yields
$$\Flip(\Der(\Cat))=\R\left(\frac{1-\sqrt{1-4t}}{2t}, \frac{t}{\sqrt{1-4t}}\right) = 
\begin{bmatrix}
1 & 0 & 0 & 0 & 0 & 0 & \cdots\\
1 & 1 & 0 & 0 & 0 & 0 & \cdots\\
2 & 3 & 1 & 0 & 0 & 0 & \cdots\\
5 & 10 & 5 & 1 & 0 & 0 & \cdots\\
14 & 35 & 22 & 7 & 1 & 0 & \cdots\\
42 & 126 & 93 & 38 & 9 & 1 & \cdots\\
\vdots & \vdots & \vdots & \vdots & \vdots & \vdots & \ddots
\end{bmatrix}.
$$ 
The row sums of $\Flip(\Der(\Cat))$ are given by the generating function
$$\frac{1 - 4t - \sqrt{1 - 4 t}}{2t (t-\sqrt{1 - 4 t})}
=\frac{1-5t+4t^2 -(1-5t)\sqrt{1-4t}}{2t(1-4t-t^2)}.$$

This is the generating function for \cite[A026737]{OEIS}, which is the number of permutations avoiding the patterns $\{3241,3421,4321\}$.


Noting that 
$$\frac{1 - 4t - \sqrt{1 - 4 t}}{2t (t-\sqrt{1 - 4 t})} = (2-C(t))\cdot \frac{1}{ \sqrt{1 - 4 t}-t}, $$
we can use Theorem~\ref{thm:latticePaths} to obtain the following identity.

\begin{identity}
\label{id:avoiding}
Let $b_n$ denote the number of permutations of length $n+1$ avoiding the patterns $\{3241,3421,4321\}$. Then
$$b_n = 2 \sum_{k=0}^n 4^{n-k}\binom{n-\frac{k+1}{2}}{n-k} - \sum_{i=0}^n \left(C_{n-i}\cdot \sum_{j=0}^{i} 4^{i-j}\binom{i-\frac{j+1}{2}}{i-j}\right).$$
\end{identity}

\subsection{The Shapiro Array} \label{subsec:shapiro}
The Shapiro array is the Riordan array $\Sha$ defined in \cite{MERLINI2017160} as follows
\begin{equation*}\label{eqn:ShapiroArray} \Sha:=\R(C(t), t C^2(t))=\R(C(t), C(t)-1) = \begin{bmatrix}
1 & 0 & 0 & 0 & 0 &  \cdots\\
1 & 1 & 0 & 0 & 0 &  \cdots\\
2 & 3 & 1 & 0 & 0 &  \cdots\\
5 & 9 & 5 & 1 & 0 &  \cdots\\
14 & 28 & 20 & 7 & 1 &  \cdots\\
\vdots & \vdots & \vdots & \vdots & \vdots & \ddots
\end{bmatrix}.
\end{equation*}

Its inverse $\Sha^{-1}$ is \cite[Theorem~3.1]{MERLINI2017160}
\begin{equation*}\label{eqn:ShapiroArrayInv}
\Sha^{-1}=\R\left(\frac{1}{1+t},\frac{t}{(1+t)^2}\right).
\end{equation*}

Now  $\Sha^{-1}_{n,k}=(-1)^{n-k}\binom{n+k}{n-k}$ \cite[Eqn~3.8]{MERLINI2017160} and hence 
    \[\Sha_{n,k}=\frac{2k+1}{n+k+1}\binom{2n}{n-k}. \]


$\Sha \cdot (1,1,2,5,14,\ldots)^t$ is known (see \cite[A007852]{OEIS}) to give the sequence $(a_n)_{n\geq 0} := (1, 2, 7, 29,$ $131,$ $625,$ $3099,\ldots)$, where $a_n$ is the number of antichains in a rooted plane tree on $n$ nodes.

We obtain the following known expression for $a_n$.
\begin{identity}
\label{id:antichains}
$$a_n = \sum_{k\geq 0} \frac{2k+1}{n+k+1}\binom{2n}{n-k} \cdot C_k = \sum_{k\geq 0} \frac{1}{2n+1}\binom{2n+1}{n-k} \binom{2k+1}{k}$$
\end{identity}

Let us now consider $$\Der(\Sha) = \R(C'(t),tC(t)) = \begin{bmatrix}
1 & 0 & 0 & 0 & 0 &  \cdots\\
4 & 1 & 0 & 0 & 0 &  \cdots\\
15 & 5 & 1 & 0 & 0 &  \cdots\\
56 & 21 & 6 & 1 & 0 &  \cdots\\
210 & 84 & 28 & 7 & 1 & \cdots \\
\vdots & \vdots & \vdots & \vdots & \vdots & \ddots
\end{bmatrix}. $$

\begin{theorem}\label{thm:nk-entriesDerSha}
For the Shapiro array $\Sha=\R(C(t), tC^2(t))$:
\begin{equation*} \Der(\Sha)=\R(C'(t),  tC(t)) \text{ and }
(\Der(\Sha))^{-1}=\R((1-t)^2(1-2t), t-t^2), \end{equation*}
and hence 
\[\Der(\Sha)_{n,k}=\binom{2n-k+2}{n-k}.\]
\end{theorem}
\begin{proof}
The proof follows along the lines of Theorem~\ref{thm:nk-entriesDerCat}.
The derivation for the $(n,k)$th entry of $\Der(\Sha)$ is as follows. Let $\Der(\Sha)=\R(d^*(t), h^*(t))$ where  
\[d^*(t)=C'(t)=\frac{1-2t-\sqrt{1-4t}}{2t^2\sqrt{1-4t}},\  h^*(t)=t\,C(t).\]  Now we compute \[\frac{1}{d^*(t-t^2)}=\frac{1}{C'(t-t^2)}=(1-t)^2(1-2t),\] and the claims follow as before from Theorem~\ref{thm:InvRiordanElements}. Once more,  
 Equation~\eqref{eqn:DerSh-cols1}  follows  using Definition~\ref{def:riordan}.
\end{proof}

This implies that, for each fixed $k\ge 0$, 
\begin{equation}\label{eqn:DerSh-cols1}
C^{k+2}(t)=2\sqrt{1-4t}\sum_{n\ge k} \binom{2n-k+2}{n+2} t^{n-k+1},
\end{equation}
or equivalently
\begin{equation}\label{eqn:DerSh-cols2}
C^{k+2}(t)\sum_{n\ge 0}\binom{2n}{n}t^n=2\sum_{n\ge k}  \binom{2n-k+2}{n+2} t^{n-k+1}.
\end{equation}

From \eqref{eqn:DerSh-cols1} and the power series for $C^k(t)$ we then obtain the following. 

\begin{identity}
\label{id:weakCompCatalan2}
\[\sum_I \binom{2i_0}{i_0}C_{i_1} C_{i_2}\ldots C_{i_{k+2}}=2\binom{2n-k+2}{n+2},\]
where the sum on the left runs over all weak compositions $I=(i_0,i_1,\ldots i_{k+2})$ of $n-k+1$, for fixed $k\ge 0.$
\end{identity}

Using Theorem~\ref{thm:FTRA}, we see that the product $\Der(\Sha)\cdot (2^i)^t_{i\geq 0}$ gives a sequence whose generating function is 
$$C'(t)\dfrac{1}{1-(1-\sqrt{1-4t})}=\dfrac{1-2t-\sqrt{1-4t}}{2 t^2(1-4t)}=\dfrac{C(t)-1}{t(1-4t)}.$$ 
Note this is the generating function for the convolution of the Catalan numbers shifted by one, and the powers of $4$, i.e., for the sequence whose $n$th term is 
$\sum_{k=1}^n C_{k+1}4^{n-k}.$
Using Theorem~\ref{thm:nk-entriesDerSha} we obtain the identity 
\begin{identity}
\label{id:derSH}
\[\sum_{k\ge 0} 2^k\Der(\Sha)_{n,k}=4^{n+1} - \binom{2n+3}{n+1} \]
or 
\[\sum_{k\ge 0} 2^k\binom{2n-k+2}{n-k}=4^{n+1} - \binom{2n+3}{n+1} \]
\end{identity}
\begin{proof} (Combinatorial) 
Let us rewrite the identity in the form
\[2\sum\limits_{k = 1}^{n} 2^{k-1}\binom{2n-k+2}{n+2}=2^{2n+2} - \left(\binom{2n+2}{n} + \binom{2n+2}{n+1} + \binom{2n+2}{n+2}\right). \]
The right-hand side counts all of the subsets of $[2n+2]$, which do not have $n,n+1$ or $n+2$ elements. Denote this set of subsets by $A$. We will show that the subsets in $A$ are counted by the left-hand side. First, note that for each subset in $A$, either $n+3$ or more of the all $2n+2$ elements are selected or $n+3$ or more of these elements are not selected. The number of subsets in these two groups is equal, so it suffices to count the subsets $Q\in A$ with at least $n+3$ elements and multiply their number by $2$. 

Let $k$ be the largest number, such that $[2n+3-k]$ contains $n+3$ of the elements in $Q$. The possible values of $k$ are between $1$ and $n$. Obviously, the element $2n+3-k \in Q$ since $k$ is the largest with this property. Thus, to determine $Q$, we shall select the $n+2$ elements in $[2n+2-k]$, which are elements of $Q$. This can be done in $\binom{2n-k+2}{n+2}$ ways. In addition, $Q$ can contain any of the $2^{k-1}$ subsets of $\{2n+4-k, 2n+5-k, \ldots , 2n+2\}$.
\end{proof}
Note that the right-hand side of Identity~\ref{id:derSH} equals the total number of inversions in all $321$-avoiding permutations of length $n$, as well as the sum of the areas of all Dyck paths of semilength $n$ \cite{cheng2007area}.


\section{Future directions}
In this article, we introduced and explored sums, derivatives, and flips of Riordan arrays. 
Using these operations, we obtained several combinatorial identities. 
We only studied the effects of these operations on a small subset of Riordan arrays, thereby merely scratching the surface of possible identities obtainable using these operations.
Furthermore, we considered only ordinary generating functions. 
A natural future direction is therefore to extend our results to the setting of exponential generating functions.
It would also be interesting to know if the second order recurrence of Theorem~\ref{thm:Eric} gives rise to combinatorially meaningful identities.

\section{Acknowledgements}

This project was carried out as part of the 2021 Graduate Research Workshop in Combinatorics.  We thank the organizers for their generous support of a valuable collaborative research experience.  


Caroline Bang thanks the National Science Foundation for support through grant DMS-1839918.

\bibliographystyle{plain}


\end{document}